\documentclass{article}

\usepackage{amssymb, theorem, amsmath, rotate}
\usepackage[all,2cell,dvips]{xy} \UseAllTwocells \SilentMatrices

\title{An elliptic K3 surface associated to Heron triangles}
\author{Ronald van Luijk}

\newcommand\luijklist{\begin{itemize}
                      \setlength{\itemsep}{-1mm}
                     }
\newcommand\eindluijklist{\end{itemize}}
\newcommand\luijkreflist{\begin{itemize}
                      \setlength{\itemsep}{-1mm}
                     }
\newcommand\eindluijkreflist{\end{itemize}}

\newcounter{nootje}
\renewcommand\check[1]{}
\newcommand\donecheck[1]{}

\newcommand\checkgone[1]{}

\newcommand\eindproof{\unskip\nobreak\hfill\hbox{\quad $\square$}\par \medskip}

\makeatletter
\def\eqalign#1{\null\,\vcenter{\openup\jot\m@th
  \ialign{\strut\hfil$\displaystyle{##}$&$\displaystyle{{}##}$\hfil
      \crcr#1\crcr}}\,}
\makeatother

\newcounter{sect}
\newcommand\mysection[1]{{\addtocounter{sect}{1}} \setcounter{theorem}{0}
                         \smallskip \noindent {\large \bf \thesect.\,\, #1} 
                         \newline \smallskip}
\newcommand\mysectiononumber[1]{{\addtocounter{sect}{1}} \setcounter{theorem}{0}
                         \smallskip \noindent {\large \bf #1} 
                         \newline \smallskip}

\renewcommand\abstract[1]
    {\begin{center}
     \small
     \begin{minipage}{4in}#1\end{minipage}
     \end{center}
    }
\newcommand\myref[1]{[#1]}
\newtheorem{theorem}{Theorem}
\newtheorem{proposition}[theorem]{Proposition}
\newtheorem{definition}[theorem]{Definition}
\newtheorem{lemma}[theorem]{Lemma}
\newtheorem{corollary}[theorem]{Corollary}
\newenvironment{proof}{\noindent {\bf Proof. }}
                      {\eindproof}
\newenvironment{proofof}{\noindent {\bf Proof of }}{\eindproof}

{\theorembodyfont{\rmfamily}

\newtheorem{remark}[theorem]{Remark}}

\renewcommand\O{\mathcal{O}}

\newcommand\T{\mathcal{T}}
\newcommand\shfHom{\mathcal{H}om}
\newcommand\Normal{\mathcal{N}}
\newcommand\Z{\mathbb{Z}}
\newcommand\Q{\mathbb{Q}}
\newcommand\F{\mathbb{F}}
\newcommand\R{\mathbb{R}}
\newcommand\G{\mathbb{G}}

\newcommand\C{\mathbb{C}}
\newcommand\oY{\overline{Y}}
\newcommand\oC{\overline{C}}
\renewcommand\P{\mathbb{P}}

\newcommand\Spec{\mathop{\rm Spec} \nolimits}
\newcommand\bSpec{\mathop{\bf Spec} \nolimits}
\newcommand\disc{\mathop{\rm disc} \nolimits}
\newcommand\Proj{\mathop{\rm Proj} \nolimits}

\newcommand\Tr{\mathop{\rm Tr} \nolimits}
\newcommand\charpoly{\mathop{\rm char} \nolimits}

\newcommand\Pic{\mathop{\rm Pic} \nolimits}
\renewcommand\deg{\mathop{\rm deg} \nolimits}
\newcommand\NS{\mathop{\rm NS} \nolimits}
\newcommand\Hom{\mathop{\rm Hom} \nolimits}

\newcommand\Num{\mathop{\rm Pic}^n \nolimits}

\newcommand\tors{{\mathop{\rm tors} \nolimits}}

\newcommand\rk{\mathop{\rm rk} \nolimits}
\newcommand\Id{\mathop{\rm Id} \nolimits}

\newcommand\Aut{\mathop{\rm Aut} \nolimits\,}
\newcommand\charac{\mathop{\rm char} \nolimits\,}
\newcommand\et{\text{\rm \'et}}
\newcommand\sep{\text{\rm sep}}

\renewcommand\frak[1]{\mathfrak{#1}}

\newcommand\ro{\xi}

\newcommand\isom{\cong}

\newcommand\ra{\rightarrow}
\newcommand\M{\tilde{M}}
\newcommand\N{\tilde{N}}
\renewcommand\tilde{\widetilde}
\newcommand\f{\mathfrak{f}}
\newcommand\g{\mathfrak{g}}
\newcommand\Mi{\mathfrak{M}}
\newcommand\Ni{\mathfrak{N}}
\newcommand\X{\mathfrak{X}}
\newcommand\Y{\mathfrak{Y}}
\newcommand\Ci{\mathfrak{C}}
\newcommand\mcol{:}
\newcommand\tnum{\tau}
\newcommand\snum{\sigma}
\newcommand\svar{s}
\newcommand\im{\mathop{\rm im} \nolimits\,}

\newcommand\myrefart[6]{\item[{[#1]}] #2, #3, {\em #4}, #5, pp. #6.}
\newcommand\myrefbook[5]{\item[{[#1]}] #2, {\em #3}, #4, #5.}

\newcommand\bogomolovtschink{BT}
\newcommand\shioda{Shi}
\newcommand\manin{Man}
\newcommand\nagata{Na}
\newcommand\hag{Ha2}
\newcommand\bruce{BW}
\newcommand\duval{Du}
\newcommand\pink{Pi}
\newcommand\sull{O'Su}
\newcommand\hagtwo{Ha1}
\newcommand\silv{Si1}
\newcommand\silvtwo{Si2}
\newcommand\aas{Aa}

\newcommand\tate{Ta3}
\newcommand\neron{BLR}
\newcommand\tatetwo{Ta1}
\newcommand\nyog{NO}
\newcommand\tatethree{Ta2}
\newcommand\maz{Maz}

\newcommand\chin{Ch}
\newcommand\artin{Ar}
\newcommand\shaf{Sha}
\newcommand\lich{Lic}
\newcommand\lip{Lip}
\newcommand\kod{Ko2}
\newcommand\kodtwo{Ko1}

\newcommand\bombmum{BM}
\newcommand\egaII{EGA II}
\newcommand\egatwo{EGA IV(2)}
\newcommand\egaone{EGA IV(1)}
\newcommand\egafour{EGA IV(4)}
\newcommand\sgaone{SGA 1}
\newcommand\sgasix{SGA 6}
\newcommand\sgafourh{SGA $4\frac{1}{2}$}
\newcommand\kramluc{KL}
\newcommand\milne{Mi}

\newcommand\deligne{De}

\newcommand\titleone{Introduction}
\newcommand\titletwo{A surface associated to Heron triangles}
\newcommand\titlethree{Elliptic surfaces}
\newcommand\titlefour{Constructions of elliptic surfaces}
\newcommand\titlefive{Proof of the main theorem}
\newcommand\titlesix{The N\'eron-Severi group under good reduction}
\newcommand\titleseven{Computing the N\'eron-Severi group and the
  Mordell-Weil group}

\newcommand\titlenine{References}

\newcommand\parag{\S}

\begin{document}

\begin{center}
\huge An elliptic K3 surface associated to Heron triangles \\
\bigskip
\large Ronald van Luijk \\
Department of Mathematics 3840 \\
970 Evans Hall \\
University of California \\
Berkeley, CA 94720-3840 \\
\verb|rmluijk@math.berkeley.edu|
\end{center}

\newpage

\noindent{\bf Abstract}:
\abstract{A rational triangle is a triangle with rational sides and
  rational area. A Heron triangle is a triangle with integral sides
  and integral area.
  In this article we will show that there exist
  infinitely many rational parametrizations, in terms of $\svar$, of
  rational triangles with perimeter $2\svar(\svar+1)$ and area
  $\svar(\svar^2-1)$. As a corollary, there exist arbitrarily many Heron
  triangles with all the same area and the same perimeter.
  The proof uses an elliptic K3 surface $Y$. Its Picard number 
  is computed to be $18$ after we prove that the N\'eron-Severi group of $Y$ injects 
  naturally into the N\'eron-Severi group of the reduction of $Y$ at a prime
  of good reduction. We also give some constructions of elliptic surfaces and 
  prove that under
  mild conditions a cubic surface in $\P^3$ can be given the structure
  of an elliptic surface by cutting it with the family of hyperplanes
  through a given line $L$. Some of these constructions were already 
  known, but appear to have lacked proof in the literature until now.
}

\bigskip

\noindent {\bf Keywords}: 
\abstract{Heron triangles, elliptic surfaces, K3
surfaces, N\'eron-Severi group, Tate conjecture, Diophantine equations, 
arithmetic geometry, algebraic geometry, number theory}

\newpage

\begin{center}
\begin{tabular}{rl}
1.& \titleone \hfill \pageref{secone}\\
2.& \titletwo \hfill \pageref{sectwo}\\
3.& \titlethree \hfill \pageref{secthree}\\
4.& \titlefour \hfill \pageref{secfour}\\
5.& \titlefive \hfill \pageref{secfive}\\
6.& \titlesix \hfill \pageref{secsix}\\
7.& \titleseven \hspace{1.2cm}\hfill \pageref{secseven}\\
  & \titlenine \hfill \pageref{secnine}\\
\end{tabular}
\end{center}

\bigskip

\frenchspacing

\mysection{\titleone}\label{secone}

A {\em rational triangle} is a triangle with rational sides and area.
A {\em Heron triangle} is a triangle with integral sides and area. 
Let $\Q(\svar)$ denote the field of rational functions in $\svar$ with
coefficients in $\Q$. The main theorem of this paper states the following.

\begin{theorem}\label{maintheorem}
There exists a sequence $\{(a_n,b_n,c_n)\}_{n\geq 1}$ of triples of
elements in $\Q(\svar)$ such that 
\luijklist
\item[{\rm 1.}] for all $n \geq 1$ and all $\snum \in \R$ with $\snum>1$,
  there exists a triangle $\Delta_n(\snum)$ with sides $a_n(\snum)$,
  $b_n(\snum)$, and $c_n(\snum)$, inradius $\snum-1$, perimeter
  $2\snum(\snum+1)$, and area $\snum(\snum^2-1)$, and
\item[{\rm 2.}] for all $m,n \geq 1$ and $\snum_0,\snum_1 \in \Q$ with
  $\snum_0,\snum_1 >1$, the
  rational triangles $\Delta_m(\snum_0)$ and $\Delta_n(\snum_1)$ are
  similar if and only if $m=n$ and $\snum_0=\snum_1$. 
\eindluijklist
\end{theorem}

\begin{remark}\label{concreteparam}
The triples of the sequence mentioned in Theorem \ref{maintheorem} can be 
computed explicitly. We will see that we can take the first four to be 
\begin{equation}\label{concrete}
(a_n,b_n,c_n) = \left(\frac{\svar(\svar+1)(y_n+z_n)}{x_n+y_n+z_n}, 
\frac{\svar(\svar+1)(x_n+z_n)}{x_n+y_n+z_n}, 
\frac{\svar(\svar+1)(x_n+y_n)}{x_n+y_n+z_n}\right),
\end{equation}
with
{ \footnotesize
$$
\eqalign{
(x_1, & y_1,z_1) = \big(1 + \svar, -1 + \svar, (-1 + \svar)\svar\big), \cr 
x_2=&(-1 + \svar)(1 + 6\svar - 2\svar^2 - 2\svar^3 + \svar^4)^3, \cr
y_2=&(-1 + \svar)(-1 + 4\svar + 4\svar^2 - 4\svar^3 + \svar^4)^3, \cr
z_2=&\svar(1 + \svar)(3 + 4\svar^2 - 4\svar^3 + \svar^4)^3, \cr 
x_3=& (-1 + \svar)(1 + 2\svar + 2\svar^2 - 2\svar^3 + \svar^4)^3\cr
&(-1 - 22\svar + 66\svar^2 + 14\svar^3 - 72\svar^4 + 30\svar^5 + 6\svar^6 - 6\svar^7 + \svar^8)^3, \cr
y_3=& (1 + \svar)(-1 + 20\svar + 68\svar^2 - 84\svar^3 + 139\svar^4 +
32\svar^5 - 224\svar^6 + \cr
&64\svar^7 + 149\svar^8 - 148\svar^9 + 
60\svar^{10} - 12\svar^{11} + \svar^{12})^3, \cr
z_3=& (-1 + \svar)\svar(5 + 10\svar + 126\svar^2 + 62\svar^3 - 225\svar^4 +
52\svar^5 + 28\svar^6 + \cr
&12\svar^7 + 27\svar^8 - 62\svar^9+38\svar^{10} - 10\svar^{11} + \svar^{12})^3, \cr
x_4=&(1 + \svar)(-1 - 62\svar + 198\svar^2 + 1698\svar^3 + 7764\svar^4 - 8298\svar^5 - 10830\svar^6 + 43622\svar^7 - 15685\svar^8  \cr
   -&45356\svar^9 - 1348\svar^{10} + 75284\svar^{11} - 13088\svar^{12} - 93076\svar^{13} + 85220\svar^{14} + 12\svar^{15} -
	  49467\svar^{16} \cr
   +&  40842\svar^{17} - 16034\svar^{18} + 2282\svar^{19} + 844\svar^{20} - 546\svar^{21} + 138\svar^{22} - 18\svar^{23} +
	    \svar^{24})^3, \cr
y_4=&(-1 + \svar)(-1 + 54\svar + 550\svar^2 - 10\svar^3 + 5092\svar^4 +  16674\svar^5 + 98\svar^6 - 51662\svar^7 + 22875\svar^8 + \cr
      41&916\svar^9 - 63076\svar^{10} + 45628\svar^{11} +
13088\svar^{12} -  63644\svar^{13} + 38884\svar^{14}+ 17668\svar^{15}
- \cr
      31&195\svar^{16}+8302\svar^{17} + 8990\svar^{18} - 9554\svar^{19} + 4476\svar^{20} - 1254\svar^{21} + 218\svar^{22} - 22\svar^{23} + \svar^{24})^3, \cr
z_4=&(-1 + \svar)\svar(-7 - 28\svar - 1168\svar^2 - 2588\svar^3 + 5170\svar^4 + 6940\svar^5 + 20176\svar^6 - 10628\svar^7 - \cr
          70&305\svar^8 + 46664\svar^9 + 85440\svar^{10} - 107832\svar^{11} + 380\svar^{12}+ 66840\svar^{13}- 46848\svar^{14}+
	13656\svar^{15} - \cr
          14&65\svar^{16} - 2796\svar^{17}+ 5712\svar^{18}-
	  5228\svar^{19} + 2738\svar^{20}- 884\svar^{21} +
	  176\svar^{22}- 20\svar^{23} 	 +\svar^{24})^3.\cr
}
$$
}
\end{remark}

Multiplying these four triples by a common denominator and
substituting only integral $\snum$, we obtain an infinite parametrized
family of quadruples of pairwise nonsimilar Heron triangles, all with
the same area and the same perimeter. 
For any positive integer $N$ we can do the same to $N$ triples of the sequence. 
We find that Theorem \ref{maintheorem} implies the following corollary. 

\begin{corollary}
For every positive integer $N$ there exists an infinite
family, parametrized by $\svar\in \Z_{>0}$, of $N$-tuples of pairwise nonsimilar
Heron triangles, all $N$ with the same area $A(\svar)$ and the same
perimeter $p(\svar)$, such that
for any two different $\svar$ and $\svar '$ the corresponding ratios
$A(\svar)/p(\svar)^2$ and $A(\svar')/p(\svar')^2$ are different.

\end{corollary}

This corollary generalizes a theorem of Mohammed Aassila 
\myref{\aas}, and Alpar-Vajk Kramer and Florian Luca \myref{\kramluc}.
Their papers give identical parametrizations to prove
the existence of an infinite parametrized family of 
{\em pairs} of Heron triangles with the same area and perimeter.
The corollary also answers the question, posed by Alpar-Vajk Kramer and 
Florian Luca and later by Richard Guy, whether triples of Heron triangles 
with the same area and perimeter exist, or even $N$-tuples with $N>3$. 
Shortly after Richard Guy had posed this question, Randall Rathbun found with a
computer search a set of $8$ Heron triangles with the same area and
perimeter. Later he found the smallest $9$-tuple. Using our methods, we can find 
an $N$-tuple for any given positive integer $N$. For example, the following
table shows $20$ values of $a,b$, and $c$ such that the triangle with sides $a$, $b$, and
$c$ has perimeter $p$ and area $A$ as given.

\begin{center}
\small
\begin{tabular}{|r|r|r|}
\hline
\multicolumn{1}{|c}{$a$}&
\multicolumn{1}{|c}{$b$}&
\multicolumn{1}{|c|}{$c$} \cr
\hline
1154397878350700583600& 2324466316136026062000& 2632653985016982326400\cr
1096939160423742636000& 2485350726331508315280& 2529228292748458020720\cr
1353301222256224441200& 2044007602377661720800& 2714209354869822810000\cr
    1326882629217053462400& 2076293397636039582000& 2708342152650615927600\cr
    1175291957596867110000& 2287901677455234640800& 2648324544451607221200\cr
    1392068029775844821400& 1997996327914674087000& 2721453821813190063600\cr
    1664717974861560418800& 1703885276761144351875& 2742914927881004201325\cr
    1159621398162242215200& 2314969007387768550000& 2636927773953698206800\cr
    1582886815525601586000& 1787918651729320350240& 2740712712248787035760\cr
    1363338670812365847600& 2031949206689694692400& 2716230302001648432000\cr
    1629738181200989059200& 1739432097243363322800& 2742347901059356590000\cr
    1958819929328111850000& 1426020908550865426800& 2726677341624731695200\cr
    2256059203526140412400& 1195069414854334519500& 2660389561123234040100\cr
    2227944754401017652000& 1213597769548172408400& 2669975655554518911600\cr
    2005582596002614412784& 1385590865209533198216& 2720344718291561361000\cr
    2462169105650632177800& 1100472310428896790000& 2548876763424180004200\cr
    2198208931289532607600& 1234160196742812482000& 2679149051471363882400\cr
    2440795514101169425200& 1105486738297174396800& 2565235927105365150000\cr
    2469616851505228370400& 1099107024377149242000& 2542794303621331359600\cr
    2623055767363274578335& 1143817472264343917040& 2344644939876090476625\cr
  \hline
  \multicolumn{3}{|l|}{$p=a+b+c=6111518179503708972000$} \cr
  \multicolumn{3}{|l|}{$A=1340792724147847711994993266314426038400000$} \cr
  \hline
\end{tabular}
\smallskip

Table 1
\end{center}

We will exhibit a bijection between the set of triples $(a,b,c)$ of sides of
(rational) triangles up to scaling and 
a subset of the set of (rational) points on a certain algebraic surface.
We will
prove Theorem \ref{maintheorem} by finding infinitely many suitable curves on this
surface. It is an {\em elliptic surface} in 
the sense of Shioda \myref{\shioda}. This will be deduced from a
generalization of the following proposition.

\begin{proposition}\label{ratsurfez}
Let $X/\C$ be a nonsingular projective surface of degree $3$ in $\P^3$
and let 
$L\subset \P^3$ be a line intersecting $X$ in three different points. Let
$\tilde{X}$ be the blow-up of $X$ in these three points. Identify
$\P^1$ with the family of planes through $L$. Then the rational map 
$X \dashrightarrow \P^1$ that sends every point of $X$ to the
plane it lies in, induces a morphism $\tilde{f} \colon \tilde{X}
\rightarrow \P^1$ that has a section. Together with any of its sections, 
$\tilde{f}$ makes 
$\tilde{X}$ into an elliptic surface over $\P^1$. 
\end{proposition}

In section 2 we will describe the surface $X$ that is used to prove
the main theorem. In section 3 we introduce the notion of elliptic
surface and state some of their properties. Two constructions of
elliptic surfaces are described in section 4. These constructions 
were already known, but the proof that they actually yield elliptic
surfaces appears to lack in the literature. Section 4 therefore
contains detailed proofs of these technical facts. One of them is a
generalization of Proposition \ref{ratsurfez}. It is used to give some
blow-up $\tilde{X}$ of $X$ the structure of an elliptic surface over
$\P^1$ in section 5. In that same section we prove 
Theorem \ref {maintheorem} by using an
elliptic K3 surface $Y \rightarrow C$, obtained from $\tilde{X}
\rightarrow \P^1$ by a base change $C \rightarrow \P^1$. 

The relation between the geometry and the arithmetic of K3 surfaces 
in general is not yet clear at all, see \myref{\bogomolovtschink}. 
The last two sections are therefore dedicated to a deeper analysis of the
geometry of the K3 surface $Y$.
They are not needed for the proof of the main theorem and serve their
own interest.
Section 6 describes the behavior of the N\'eron-Severi group of a surface
under good reduction. This again was already known, but until
now lacked proof in the literature. It is used in section 7 to
determine the full N\'eron-Severi group of $Y$ and the Mordell-Weil group
of the generic fiber of $Y \rightarrow C$. 

The author would like to thank Bjorn Poonen,
Arthur Ogus, Robin Hartshorne, Tom Graber, Bas Edixhoven, Jasper Scholten,
and especially Hendrik Lenstra for very helpful discussions.
\bigskip

\mysection{\titletwo}\label{sectwo}

For a triangle with sides $a$, $b$, and $c$, let $r$, $p$, and $A$
denote its inradius, perimeter, and area respectively. The line
segments from the vertices of the triangle to the midpoint of the
incircle divide the triangle in three 
smaller triangles of areas $ar/2$, $br/2$, and
$cr/2$. Adding these we find $A=rp/2$. Set $x=p/2-a$, 
$y=p/2-b$, and $z=p/2-c$. Then we get $p=2(x+y+z)$, so $A=r(x+y+z)$. 
Heron's formula $A^2=(x+y+z)xyz$ then yields $r^2(x+y+z)=xyz$. Therefore,   
the point $[r \mcol x \mcol y\mcol z]\in \P^3$ lies on
the surface $X \subset \P^3_\Q$ given by $r^2(x+y+z) = xyz$. Conversely, if  
$[1 \mcol x \mcol y\mcol z]$ lies on $X$, with $x, y, z >0$, then
the triangle with sides $a=y+z$, $b=x+z$, and $c=x+y$ has inradius $1$. 
Thus we get a bijection between the set of triples $(a,b,c)$ of sides
of triangles up to scaling and the set of real points $[r \mcol x \mcol y
  \mcol z]$ on $X$ with positive ratios $x/r$, $y/r$, and $z/r$.
Let $G \subset \Aut X$ 
denote the group of automorphisms of $X$ induced by the permutations
of the coordinates $x$, $y$, and $z$. Let $f \colon X \dashrightarrow
\P^1$ be the rational map given by  
$f \colon [r:x:y:z] \mapsto [r : x+y+z]$. Note that if we let $G$ act
trivially on $\P^1$, then $f$ commutes with the action of $G$. 

\begin{lemma}\label{TrianglesToPoints}
For $i=1,2$, let $\Delta_i$ denote a triangle, let $a_i$, $b_i$,
and $c_i$ denote the sides of $\Delta_i$,
and let $P_i$ be the point on $X$ corresponding to the
equivalence class (under scaling) of the triple $(a_i,b_i,c_i)$. Then
$\Delta_1$ and $\Delta_2$ are similar if and only if $P_1$ and $P_2$
are in the same orbit under $G$. 
Up to scaling, $\Delta_1$ and $\Delta_2$ have the same inradius and
perimeter if and only if $P_1$ and $P_2$ map to the same point under
$f$. 
\end{lemma}
\begin{proof}
This is obvious.
\end{proof}

To set our strategy for proving Theorem \ref{maintheorem}, note that
it asserts that for fixed $\snum$, the
infinitely many pairwise nonsimilar triangles $\Delta_n(\snum)$, with
$n\geq 1$, all have the same perimeter $2\snum(\snum+1)$ and
inradius $\snum-1$.  
By Lemma \ref{TrianglesToPoints} this is equivalent to the statement that
the infinitely many points corresponding to the triples
$(a_n(\snum),b_n(\snum),c_n(\snum))$ all map under $f$ to
$[\snum-1 : \snum(\snum+1)]$, and that they are all in different
orbits under $G$. To prove Theorem \ref{maintheorem}, we will find 
a suitable infinite collection of 
curves on $X$, mapping surjectively to $\P^1$ under
$f$. Those maps will not be surjective on rational points, but 
for rational $\snum$ each of
these curves will intersect $f^{-1}([\snum-1 : \snum(\snum+1)])$ 
in a rational point. 

\begin{remark}\label{XisRat}
Since the equation $r^2(x+y+z) = xyz$ is linear in $x$, we find that
$X$ is rational. A parametrization is given by the birational
equivalence $\P^2 \dashrightarrow X$, given by 
\begin{equation}\label{Xrational}
\eqalign{
[r:x:y:z] &= [vw(u-v):v(uv+w^2):w^2(u-v):uv(u-v)], \quad \mbox{or}\cr
[u:v:w] &= [yz:r^2:yr]. \cr
}
\end{equation}
\end{remark}

%
%
\bigskip

\mysection{\titlethree}\label{secthree}

In this section, $k$ will denote an algebraically closed field. All
varieties, unless stated otherwise, are $k$-varieties. A variety $V$
over a field $l$ is called smooth if the morphism $V \rightarrow \Spec l$
is smooth. 
We will start with the definition of a lattice. Note that for any abelian
groups $A$ and $G$, a symmetric bilinear pairing $A \times A
\rightarrow G$ is called nondegenerate if the induced homomorphism 
$A \rightarrow \Hom(A,G)$ is injective. We do not require a
lattice to be definite, only nondegenerate.

\begin{definition}
A {\em lattice} is a free $\Z$-module $L$ of finite rank, endowed with a
symmetric, bilinear, nondegenerate pairing $\langle
\underline{\,\,\,\,}\,,\underline{\,\,\,\,} \rangle \colon L \times L
\rightarrow 
\Q$. An {\em integral lattice} is a lattice with a $\Z$-valued pairing. 
A {\em sublattice} of $L$ is a submodule $L'$ of $L$, such that the induced
bilinear pairing on $L'$ is nondegenerate. The positive- or
negative-definiteness or signature of a lattice is defined to be that 
of the vector space $L_{\Q}$, together with the induced pairing.
\end{definition}

\begin{definition}
The \/ {\em Gram matrix} of a lattice $L$ with respect to a given basis
$x=(x_1,\ldots, x_n)$ is $I_x = (\langle x_i,x_j \rangle)_{i,j}$.
The {\em discriminant} of $L$ is defined by
$\disc L = \det I_x$ for any basis $x$ of
$L$. A lattice $L$ is called {\em unimodular} if $\disc L =\pm 1$.
\end{definition}

\begin{definition}
A {\em fibration} of a variety $Y$ {\em over} a 
regular integral curve $Z$ over $k$ 
is a dominating morphism $g \colon Y \rightarrow Z$.
\end{definition}

\begin{remark}\label{flatfibration}
If $Y$ is integral in the definition above, then $g$ is flat, see
\myref{\hag}, Prop. III.9.7. If also the characteristic of $k$ equals $0$
and the singular locus of
$Y$ is contained in finitely many fibers, then almost
all fibers are nonsingular, see \myref{\hag}, Thm. III.10.7.
If $Y$ is projective, then $g$ is surjective, as projective morphisms
are closed.
\end{remark}

\begin{definition}
A fibration of a smooth, projective, irreducible surface $Y$ over 
a smooth, projective, irreducible curve $Z$ is called 
{\em relatively minimal} if for every 
fibration of a smooth, projective, irreducible surface $Y'$ over $Z$, every
$Z$-birational morphism $Y \rightarrow Y'$ is necessarily an isomorphism.
\end{definition}

\begin{theorem}\label{blowdownable}
Let $Y$ be a smooth, projective, irreducible surface, $Z$ a smooth,
projective, irreducible curve, and let $g \colon Y \rightarrow Z$ be a
fibration  
such that no fiber contains an exceptional prime divisor $E$, i.e., 
a prime divisor with self-intersection number $E^2=-1$ and
$H^1(E,\O_E)=0$. Then $g$ is a relatively minimal fibration. 
\end{theorem}
\begin{proof}
This is a direct corollary of the Castelnuovo Criterion
(\myref{\chin}, Thm. 3.1)  
and the Minimal Models Theorem (\myref{\chin}, Thm. 1.2). See also Lichtenbaum
\myref{\lich} and Shafarevich \myref{\shaf}.
\end{proof}

%

\begin{lemma}\label{connfibers}
Let $g\colon X \rightarrow Y$ be a projective morphism of noetherian
schemes. Assume that $X$ is integral and that $g$ has a section. Then
there is an isomorphism $g_* \O_X \isom \O_Y$ if and only if for every
$y \in Y$ the fiber $X_y$ is connected.
\end{lemma}
\begin{proof}
Set $Y' = \bSpec g_* \O_X$.
By Stein factorization (see \myref{\hag}, Cor. III.11.5) the morphism 
$g$ factors naturally as $g = h \circ f$,
where $f \colon X \rightarrow Y'$ is projective with connected 
fibers and $h \colon
Y' \rightarrow Y$ is finite. If we have $g_* \O_X \isom \O_Y$, then
$h$ is an isomorphism, so $g$ has connected fibers. Conversely, suppose $g$
has connected fibers. As $f$ is projective, it is closed. If $f$ were
not surjective, then there would be a nonempty open affine 
$V \subset Y'$ with $f^{-1}(V)=\emptyset$. This implies
$(f_*\O_X)(V)=0$, contradicting the equality $f_* \O_X
= \O_{Y'}$. We conclude that $f$ is surjective, so $h$ also has
connected fibers. As $h$ is finite, its fibers are also totally
disconnected (see \myref{\hag}, exc. II.3.5), so $h$ is injective on
topological spaces. Let $\varphi \colon Y \rightarrow X$
be a section of $g$. Then $\psi=f \circ \varphi$ is a section of
$h$. Every injective continuous map between topological spaces that
has a continuous section is a homeomorphism, so $h$ is
a homeomorphism. Therefore, to prove that $h$ is
an isomorphism, it suffices to show this locally, so we may
assume $Y'=\Spec B$ and $Y=\Spec A$. The composition $\psi^{\#} \circ
h^{\#} \colon A \rightarrow B \rightarrow A$ is the identity, so
$\psi^{\#}$ is surjective. As $X$ is integral, so is $Y'$. Hence, the
ideal $(0) \subset B$ is prime. 
Since $\psi$ is surjective, there is a prime ideal $\frak{p}
\subset A$ such that $(0)=\psi(\frak{p}) = (\psi^{\#})^{-1} \frak{p}$,
so $\psi^{\#}$
is injective. We find that $\psi^{\#}$ is an isomorphism. Hence,
so are $\psi$ and $h$, so there is an isomorphism $g_* \O_X \isom \O_Y$.
\end{proof}

\begin{definition}
A fibration is called \/ {\em elliptic} if all but finitely many fibers
are curves of genus $1$.
\end{definition}

\begin{theorem}{\label{canonicalonminimalfibration}}
Let $C$ be a smooth, irreducible, projective curve of genus $p(C)$ over
an algebraically closed 
field $k$. Let $S$ be a smooth, irreducible, projective surface over
$k$ with Euler characteristic $\chi=\chi(\O_S)$ and
let $g \colon S \rightarrow C$ be an elliptic
fibration that has a section. Then the following are equivalent.
\luijklist
\item[{\rm (i)}] The morphism $g$ is a relatively minimal fibration, 
\item[{\rm (ii)}] Any canonical divisor $K_S$ on $S$ is algebraically equivalent to 
$(2p(C)-2+\chi)F$, where $F$ is any fiber of $g$, 
\item[{\rm (iii)}] We have $K_S^2=0$.
\eindluijklist
\end{theorem}
\begin{proof}
(i) $\Rightarrow$ (ii). By Remark \ref{flatfibration} the morphism $g$
  is flat, so by the principle of connectedness, {\em all} fibers are
  connected, see \myref{\hag}, exc. III.11.4. From Lemma
  \ref{connfibers} we find that $g_* \O_S \isom \O_C$. 
  Under that assumption, an explicit expression for $K_S$ is given in 
  \myref{\kodtwo}, \parag{} 12, for base fields that can be embedded in $\C$, 
  and in \myref{\bombmum}, \parag{} 1, for characteristic $p$.
Since $g$ has a section, say $\O$, every fiber of $g$ will have
intersection multiplicity $1$ with the horizontal divisor $\O(C)$, so
there are no multiple fibers. In that case, the expression mentioned
above implies that $K_S$ is
algebraically equivalent to $(2p(C)-2+\chi)F$ for any fiber $F$. 

(ii) $\Rightarrow$ (iii).
Since $F$ is algebraically equivalent to any other fiber $F'$, 
it is also numerically equivalent to any other fiber $F'$.
Thus we get $F^2=F \cdot F'=0$, so $K_S^2=0$.

(iii) $\Rightarrow$ (i). 
If $g$ were not relatively minimal, then the Minimal Models Theorem
(see \myref{\chin}, Thm. 1.2) tells us that there would be a
relatively minimal fibration $g' \colon S' \rightarrow C$ of a smooth,
irreducible, projective surface $S'$ and a
$C$-morphism $\gamma \colon S\rightarrow S'$ which consists of a nonempty
sequence of blow-ups of points. Then $g'$ is an elliptic fibration
as well. The composition $\gamma \circ \O$ is a section of $g'$. 
By the proven implication (i) $\Rightarrow$ (iii), we find that
$K_{S'}^2=0$. This implies $K_{S}^2 <0$, because
for any blow-up $Z \rightarrow Z'$ of a nonsingular projective surface
$Z'$ in a point $P$, we have $K_{Z}^2=K_{Z'}^2-1$, see \myref{\hag},
Prop. V.3.3. From this contradiction, we conlude that $g$ is
relatively minimal.
\end{proof}

The following definition states that 
if the fibration $g$ as described in Theorem
\ref{canonicalonminimalfibration} is not smooth, then we call the
quadruple $(S,C,g,\O)$ an {\em elliptic surface}. Note that throughout 
this section $k$ is assumed to be algebraically closed.

\begin{definition}\label{ellsurf}
Let $C$ be a smooth, irreducible, projective curve over $k$. 
An {\em elliptic surface over} $C$ is a smooth, irreducible, projective 
surface $S$ over $k$ together with a relatively minimal elliptic fibration 
$g \colon S \rightarrow C$ that is {\bf not smooth}, and a 
section $\O \colon C \rightarrow S$ of $g$. 
\end{definition}

\begin{remark}\label{smooth}
In order to rephrase what it means for $g$ not to be smooth, note that by
\myref{\egatwo}, D\'ef. 6.8.1, a morphism of schemes $g
\colon X \rightarrow Y$ is smooth if and only if $g$ is flat, $g$ is
locally of 
finite presentation, and for all $y \in Y$ the fiber $X_y = X \times_Y
\Spec k(y)$ over the residue field $k(y)$ is geometrically regular.
See also \myref{\hag}, Thm. III.10.2.

In the case that $g$ is a fibration of an integral variety $X$ over a 
smooth, irreducible, projective
curve over an algebraically closed field $k$, it follows from Remark
\ref{flatfibration} that $g$ is flat. As $X$ is noetherian and of 
finite type over $k$, it also follows that $g$ is locally of finite
presentation.  Hence $g$ not being smooth is then equivalent to the
existence of a singular fiber. 
\end{remark}

For the rest of this section,
let $S$ be an elliptic surface over a smooth, irreducible, projective curve
$C$ over $k$, fibered by $g\colon S \ra C$ 
with a section $\O$. Let $K=k(C)$ denote the function field of $C$ 
and let $\eta \colon \Spec K\rightarrow C$ be its generic point. Then the generic fiber 
$E=S \times_C \Spec K$ of $g$ is a curve over $K$ of genus $1$. 
Let $\xi$ denote the natural map $E\rightarrow S$. 

$$
\xymatrix{ E \ar[r]^{\xi} \ar[d] & S \ar[d]^g \\
         \Spec K  \ar[r]_{\eta} &C}
$$

\begin{lemma}\label{sections}
The identity on $\Spec K$ and 
the composition of any section of $g$ with $\eta$, together 
induce a section of $E \rightarrow \Spec K$. 
This correspondence induces a bijection 
between the set $E(K)$ of sections of $E\ra \Spec K$
and the set $S(C)$ of sections of $g$. 
\end{lemma}
\begin{proof}
The first sentence follows from the universal property of fibered products.
As $\Spec K$ is dense in $C$,  
any section of $E \rightarrow \Spec K$, composed with $\xi$, 
induces a section of $g$ on an open 
part of $C$. As $C$ is a smooth curve and $S$ is projective, this extends 
uniquely to a section of $g$, see \myref{\hag}, Prop. I.6.8.
This map is clearly the inverse of the map described in the lemma.
\end{proof}

Whenever we implicitly identify the two sets $E(K)$ and $S(C)$, we
will do this using the bijection of Lemma \ref{sections}.
The section $\O$ of $g$ corresponds to a point on $E$, giving $E$ 
the structure of an elliptic curve. This puts a group structure on
$E(K)$, which carries over to $S(C)$, see \myref{\silv}, Prop. III.3.4.  

Recall that for any proper scheme $Y$ over an algebraically closed field, 
the N\'eron-Severi group $\NS(Y)$ of $Y$ is the quotient of $\Pic Y$ by
the group $\Pic^0 Y$ consisting of all divisor classes algebraically 
equivalent to $0$. For a precise definition of algebraic equivalence,
see \myref{\hag}, exc.~V.1.7, which is stated only for smooth surfaces, but
holds in any dimension, see \myref{\sgasix}, Exp. XIII, p. 644, 4.4. 
\donecheck{Find reference for alg. equiv. over nonalg. closed field, and 
nonproper, then rewrite this parag.}
We will write $D \sim D'$ and $D \approx D'$ 
to indicate that two divisors $D$ and $D'$ are linearly and algebraicallly 
equivalent respectively.
Algebraic equivalence implies numerical equivalence, see
\myref{\sgasix}, Exp. X, p. 537, D\'ef. 2.4.1, and p. 546, Cor. 4.5.3.
If $Y$ is proper, then $\NS(Y)$ is a
finitely generated, abelian group, see \myref{\hag}, exc.~V.1.7--8,
or \myref{\milne}, Thm. V.3.25 for surfaces, or \myref{\sgasix},
Exp. XIII, Thm. 5.1 in general.
Its rank $\rho=\rk \NS (Y)$ is called the Picard number of $Y$.

\begin{proposition}
On $S$, algebraic equivalence coincides with numerical equivalence. The group
$\NS(S)$ is torsion-free. The intersection pairing induces a
nondegenerate bilinear pairing on $\NS(S)$, making it into a lattice
of signature $(1,\rho-1)$. 
\end{proposition}
\begin{proof}
The first statement is proven by Shioda in \myref{\shioda}, Thm.~3.1.
It follows immediately that $\NS(S)$ has no torsion and that the
bilinear intersection pairing is nondegenerate on $\NS(S)$, see
\myref{\shioda}, Thm.~2.1 or \myref{\hag}, example V.1.9.1. The
signature follows from the Hodge Index Theorem (\myref{\hag}, Thm.~V.1.9).
\end{proof}

\begin{lemma}\label{samepiczero}
The induced map $g^* \colon \Pic^0 C \ra \Pic^0 S$ is an isomorphism. 
\end{lemma}
\begin{proof}
See \myref{\shioda}, Thm.~4.1. 
\end{proof}

The map $\ro$ induces a homomorphism $\ro^* \colon\Pic S \ra \Pic E$. 
Since $\ro^* \circ g^* = (g\circ \ro)^*$ 
factors through $\Pic (\Spec K)=0$, we have an inclusion 
$g^*\Pic^0 C \subset \ker \ro^*$, so $\ro^*$ factors through 
$\Pic S/g^*\Pic^0 C$, which is isomorphic to $\NS(S)$ by Lemma
\ref{samepiczero}. Let 
$\varphi \colon \NS(S) \rightarrow E(K)$ denote the composition of the
induced homomorphism $\NS(S) \ra \Pic E $ with the homomorphism 
$\Sigma \colon \Pic E  \ra E(K)$, which uses the group law on $E(K)$ to 
add up all the points of a given divisor, with multiplicities. 

$$
\xymatrix{\Pic S \ar[r] \ar[d]_{\xi^*} & \NS(S) \ar[ld] \ar[d]^\varphi \\
          \Pic E \ar[r]_{\Sigma}        & E(K)}
$$

Set $T=\ker \varphi$ and for $v \in C$, let $m_v$ denote the number of 
irreducible components of the fiber of $g$ at $v$.
Let $r$ denote the rank of the Mordell-Weil group $E(K)$. Finally,
for every point $P \in E(K)$, let $(P)$ denote the prime divisor on $S$ 
that is the image of the section $C \ra S$ corresponding to $P$ by Lemma 
\ref{sections}. 

\begin{proposition}\label{Tkern}
The homomorphism $\varphi$ is surjective and maps $(P)$ to $P$.
The group $T$ is a sublattice of $\NS(S)$, generated by $(\O)$ and
the irreducible components of the singular fibers of $g$.
Its rank equals $\rk T =2+\sum_v (m_v-1)$. We have
$\rho = r+2+\sum_v (m_v-1)$.
\end{proposition}
\begin{proof}
For the first claim, see \myref{\shioda}, Lemma 5.1. For the
description of the kernel $T$ and the fact that it is a lattice, 
see \myref{\shioda}, Thm.~1.3. For its
rank, see \myref{\shioda}, Prop. 2.3. The last equality then follows from
the exact sequence
\begin{equation}\label{TNsE}
0 \longrightarrow T \longrightarrow \NS(S) 
\longrightarrow E(K) \longrightarrow 0.
\end{equation}
\end{proof}

Shioda also shows that $\NS(S)$ is the direct sum of a negative definite
lattice of rank $\rho-2$ and the unimodular
lattice $U$ of rank $2$ generated by $(\O)$ and $F$, where $F$ is any
fiber. Since $U$ is contained in $T=\ker \varphi$, it follows from
(\ref{TNsE}) 
that $E(K)$ is the quotient of a negative definite lattice by a sublattice. 
By general theory of definite lattices,
the nondegenerate pairing coming from the orthocomplement of $U$ in
$\NS(S)$ induces a nondegenerate pairing on $E(K)/E(K)_{\tors}$. 
Shioda, \myref{\shioda}, Thm.~8.6,
gives an explicit formula for the negative of this pairing, 
which is twice the canonical height pairing.  
\bigskip

\mysection{\titlefour}\label{secfour}

In this section we will prove that under mild conditions a fan
of hyperplane sections of a degree $3$ surface
in $\P^3$ gives an elliptic surface. 
This statement is well known, at least for nonsingular surfaces in
characteristic $0$, but details such as the
existence of singular fibers are often overlooked.
Also under mild conditions a base
extension of an elliptic surface gives again an elliptic surface. 
Both statements seem to lack proofs in the literature, so we include them here.

\begin{definition}
A surface $X$ over an algebraically closed field $k$ has 
a {\em rational singularity} at a point $x$ 
if there exist a surface $Y$ and a projective, birational morphism 
$f \colon Y \rightarrow X$ that is an isomorphism from
$f^{-1}(X-\{x\})$ to $X-\{x\}$ and such that we have $R^1f_* \O_Y =0$ and
$f^{-1}(U)$ is smooth over $k$ for some open neighborhood $U$ of $x$.
\end{definition}

\begin{remark}\label{artin}
Let $f \colon Y \rightarrow X$ be a resolution of a singularity at $x$
on $X$ with exceptional curve (possibly reducible) $E$. Then $x$ is a
rational singularity if and only if for every positive divisor $Z$ on
$Y$ with support in $E$ the arithmetic genus $p_a(Z)$ satisfies
$p_a(Z)\leq 0$, see \myref{\artin}, Prop. 1.
\end{remark}

\begin{proposition}\label{rationalellipticsurface}
Let $k$ be any field of characteristic not equal to $2$ or $3$,
contained in an algebraically closed field $k'$. Let $X$ be
a projective, irreducible surface in $\P^3_k$ of degree $3$, 
which is geometrically
regular outside a finite number of rational singularities. Let $L$ be
a line that intersects $X$ in three different nonsingular
points $M_1$, $M_2$, and $M_3$. Identify $\P^1$ with the
family of hyperplanes in $\P^3$ through $L$ and let $f \colon X
\dashrightarrow \P^1$ be the rational map that sends every point of
$X$ to the hyperplane it lies in. Let $\pi \colon \tilde{X}
\rightarrow X$ be a minimal desingularization of the blow-up of $X$ at 
the $M_i$. For $i=1,2,3$, let $\M_i$ denote the exceptional curve above
$M_i$ on $\tilde{X}$. Then $f \circ \pi$ extends to a morphism
$\tilde{f}$. It maps
the $\M_i$ isomorphically to $\P^1$, yielding at least three sections.
Together with any of its sections, $\tilde{f}$ makes
$\tilde{X}_{k'}$ into a rational elliptic surface over $\P^1_{k'}$. 
\end{proposition}

\begin{remark}
O'Sullivan (\myref{\sull}, Prop. VI.1.1) shows that any normal cubic
surface in $\P^3$ that is not a cone has only rational double points. 
He excludes characteristics $2$, $3$, and $5$, but describes how his 
results could be extended to any characteristic using results from Lipman 
\myref{\lip}. For a published reference, see \myref{\bruce}
(characterisitic $0$). 
\end{remark}

The proof of Proposition \ref{rationalellipticsurface} consists of
several steps. For clarity, we will prove them in separate lemmas.
Let $k,k',L,X,\tilde{X},\pi,M_i,\M_i,f,$ and $\tilde{f}$ be as in
Proposition \ref{rationalellipticsurface}. 
First we will show that $\tilde{X}$ is rational, smooth, and irreducible. 
Then we show that $\tilde{f}$ is a morphism that has a section. 
We proceed by showing that almost all fibers are nonsingular of genus
$1$. After that, we see that $\tilde{f}$ is not smooth and finally, we will
show that $\tilde{f}$ is a relatively minimal fibration. Note that if
$L$ is defined over $k$, then so is $f$. If $M_i$ is a $k$-point, then
the section $\O$ is defined over $k$ as well. All other statements are
geometric, so we may assume that $k=k'$.

\begin{lemma}\label{ratsurfone}
Under the assumptions of Proposition \ref{rationalellipticsurface} the
surface $\tilde{X}$ is rational, smooth, and irreducible.
\end{lemma}
\begin{proof}
By construction, $\tilde{X}$ is smooth. It is irreducible because $X$
is, and $\pi \colon \tilde{X} \rightarrow X$ is birational.
Obviously, to show that $\tilde{X}$ is rational, 
it suffices to show that $X$ is rational.
It is a classical result that nonsingular cubics are obtained by
blowing up $6$ points in general position in $\P^2$, whence they 
are rational. For this statement, see \myref{\hag}, \parag{} V.4, in
particular Rem. V.4.7.1. Proofs are given
in \myref{\manin}, \parag{} 24 or \myref{\nagata}, I, Thm. 8, p. 366.

For the singular case, note that $X$ is not a cone. Indeed, the
exceptional curve $E$ of the desingularization of a cone over a plane
cubic is isomorphic to that cubic, see \myref{\hag},
exc. II.5.7. Hence, it would satisfy $p(E)=1$, which contradicts
Remark \ref{artin}. As $X$ is not a cone, projection from any singular
point $x$ will give a dominant rational map from $X$ to $\P^2$. It is
birational because every line through $x$ that is not contained in $X$
intersects $X$ by Bezout's Theorem in only one more point.
%
%
\end{proof}

\begin{lemma}\label{ratsurfonehalf}
The rational map $\tilde{f}$ extends to a morphism, mapping 
$\M_i$ isomorphically to $\P^1$. 
\end{lemma}
\begin{proof}
The rational map $f$ is defined everywhere, except at the $M_i$,
whence the composition $f \circ \pi$ is well-defined outside the $\M_i$.
Any point $P$ on $\tilde{M}_i$ corresponds to a direction at $M_i$ on $X$. 
Since $L$ intersects $X$ in three different points and the total
intersection $L\cdot X$ has degree $3$ by B\'ezout's Theorem, it
follows that $L$ is not tangent to $X$, so these directions at $M_i$
are cut out by the planes through $L$. The map $f \circ \pi$ extends
to a morphism $\tilde{f}$ by sending $P \in \M_i$ to the plane that
cuts out the direction at $M_i$ that $P$ corresponds to. 
Thus, it induces an isomorphism from the $\M_i$ to $\P^1$.
%
%
\end{proof}

Note that if a hyperplane $H$ does not contain any singular points of $X$, 
then the fiber of $\tilde{f}$ above $H$ is isomorphic to $H \cap X$.
Here the missing
points $M_i$ in $f^{-1}(H)=(H \cap X)\setminus \{M_1,M_2,M_3\}$ are
filled in by the appropriate points on $\M_i$. To prove that almost
all fibers are nonsingular curves of genus $1$ we will use 
Proposition \ref{geomint}. Its proof was suggested by B. Poonen.

\begin{lemma}\label{conregint}
Any connected, regular variety is integral.
\end{lemma}
\begin{proof}
Let $Z$ be such a variety. Then $Z$ is reduced, so it
suffices to show that $Z$ is irreducible.
The minimal primes of the local ring of a point on $Z$
correspond to the components it lies on. As a regular
local ring has only one minimal prime ideal, we conclude that 
every point of $Z$ lies on exactly one component. As $Z$ is
connected, $Z$ is irreducible.
\end{proof}

\begin{proposition}\label{geomint}
Let $Y$ be a geometrically connected, regular variety over a field $l$.
If $Y$ contains a closed point of which the residue field is separable
over $l$, then $Y$ is geometrically integral.
\end{proposition}
\begin{proof}
Let $l^{\sep}$ denote the separable closure of $l$. As separable
extensions preserve regularity (see \myref{\egatwo}, Prop. 6.7.4), 
we find that $Y_{l^{\sep}}$ is regular. As it is connected as well, 
$Y_{l^{\sep}}$ is integral by Lemma \ref{conregint}, whence
irreducible. Over a separably
closed field, irreducibility implies geometric irreducibility, see
\myref{\egatwo}, Prop. 4.5.9. Therefore $Y_{\overline{l}}$ is irreducible.

Let $c$ be the closed point mentioned. 
Then the local ring $\O_{Y,c}$ is regular, with residue field
separable over $l$. From \myref{\egaone}, Thm. 19.6.4, we find that 
the ring $\O_{Y,c}$ is formally
smooth over $l$. By \myref{\egatwo}, Thm. 6.8.6, this means that $Y$
is smooth (over $l$) at $c$. As smoothness is an open condition (see
\myref{\egatwo}, Cor. 6.8.7), there is a nonempty open subset $U \subset Y$
such that $Y$ is smooth at all $x \in U$. As smoothness is a local
condition, $U$ is smooth, whence geometrically regular. 

As $Y_{\overline{l}}$ is irreducible, the subset $U_{\overline{l}}$ is
dense and also irreducible, whence connected. It is also regular, so it is
integral by Lemma \ref{conregint}. 
Therefore, $U$ is geometrically integral, which for an integral scheme
over $l$ is equivalent to the fact that its
function field is a primary and separable field extension of
$l$, see \myref{\egatwo}, Cor. 4.6.3. 
As $Y$ is integral and the function field $k(Y)$ of $Y$ is isomorphic to
the function field $k(U)$ of $U$, it follows that $Y$ is geometrically
integral as well. 
\end{proof}

\begin{lemma}\label{ratsurftwo}
Under the assumptions of Proposition \ref{rationalellipticsurface}
almost all fibers are nonsingular curves of genus $1$.
\end{lemma}
\begin{proof}
It follows from Remark \ref{flatfibration} that almost all fibers are 
nonsingular if $\charac k =0$. Suppose $\charac k = p >3$.
We will first show that the generic fiber $E = \tilde{X} \times_{\P^1}
\Spec k(t)$ above the generic point $\eta: \Spec k(t) \rightarrow
\P^1$ of $\P^1$ is regular. 
Then we will show $E$ is geometrically integral of genus $1$ 
and finally we will conclude it is smooth over $\Spec k(t)$.

Take a point $P \in E$ and let $x \in \tilde{X}$ be the image of $P$
under the projection $\varphi \colon E \rightarrow \tilde{X}$.  
On every open $U=\Spec A \subset \P^1$, the map $\eta$
is given by the localization map $\psi \colon A \hookrightarrow k(t)$. As
fibered products of affine spaces come from tensor products, which
commute with localization, the map $\varphi^{\#}\colon
\O_{\tilde{X},x} \rightarrow \O_{E,P}$ on local $A$-algebras is
induced by $\psi$. The maximal ideal of
$\O_{\tilde{X},x}$ pulls back under $\tilde{f}^{\#}|_A \colon A
\rightarrow \O_{\tilde{X},x}$ to the prime ideal of $A$ corresponding to
$\tilde{f}(x)=\im \eta$, i.e., to $(0)$. Hence, all nonzero
elements of $A$ are already invertible in $\O_{\tilde{X},x}$, so the
map $\O_{\tilde{X},x} \rightarrow \O_{E,P}$ is in fact an isomorphism.
Since $\tilde{X}$ is regular by Lemma \ref{ratsurfone}, we conclude that 
$\O_{\tilde{X},x} \isom \O_{E,P}$ is a regular local ring, so $E$ is regular.


Also, for any extension field $F$ of $k(t)$ the scheme $E
\times_{k(t)} F$ is a cubic in $\P^2_F$, so it is connected. Thus,
$E$ is geometrically connected. As in Lemma \ref{sections}, the
sections $\M_i$ determine $k(t)$-points on $E$. From Proposition
\ref{geomint} we find that $E$ is geometrically integral.
As $E$ is a regular, geometrically integral, plane cubic curve, 
it has genus $g(E)=1$. 
Here we define the genus $g(C)$ of a regular (but possibly not
smooth), projective, and geometrically
integral curve $C$ by the common value of its arithmetic genus
$p_a(C)$ and its geometric genus $p_g(C) = \dim
H^0(C,\omega_C^\circ)$, where
$\omega_C^\circ$ is the dualizing sheaf of $C$, see \myref{\hag},
III.7. \donecheck{see section 3.2 in Bjorn's notes}

Now, if $E$ were not smooth over $k(t)$, then there would be a finite 
extension $F/k(t)$ such that $E_F = E \times_{k(t)} F$ is not 
regular. Any nonregular plane cubic has genus $0$, so $g(E_F) = 0$. 
Let $K/k(t)$ be the subfield of $F$ such that $K/k(t)$ is
separable and $F/K$ is purely inseparable. Then by \myref{\egatwo},
Prop. 6.7.4, the curve $E_K= E \times_{k(t)} K$ is regular, 
so $g(E_K)=1$. By \myref{\tatetwo}, Cor. 1, the
difference $g(E_K)-g(E_F)=1$ is an integral multiple of $(p-1)/2$, so we
find $p=2$ or $p=3$. Since $p>3$, 
we conclude that $E$ is smooth over $\eta$. 
As $\tilde{f}$ is flat and projective, by \myref{\hag}, exc. III.10.2,
there is a dense open subset $U \subset \P^1$
on which $\tilde{f} \colon \tilde{f}^{-1}(U) \rightarrow U$ is
smooth. By Remark \ref{smooth} almost all fibers are then
nonsingular. As they are plane cubics, they have genus $1$. 
\end{proof}

\begin{lemma}\label{ratsurfthree}
Under the assumptions of Proposition \ref{rationalellipticsurface}
the morphism $\tilde{f}$ is not smooth.
\end{lemma}
\begin{proof}
By Remark \ref{smooth}, it suffices to prove that there exists a
singular fiber. 
As there are only finitely many singular points on $X$, for almost all
planes $H$ through $L$ the fiber $\tilde{X}_H$ is isomorphic to $X
\cap H$. As any two curves in $H \isom \P^2$ intersect, it
follows that $\tilde{X}_H$ is connected for all but finitely many $H$.
Since $\tilde{f}$ is flat (see Remark \ref{flatfibration}), it follows
from the principle of connectedness (see \myref{\hag},
exc. III.11.4) that the fiber $\tilde{X}_H$ is connected for all $H$.

If $X$ contains a singular point, then the fiber $\tilde{X}_H$ of
$\tilde{f}$ above the plane $H$ that it lies in contains an
exceptional curve, so it is reducible and connected.
From Lemma \ref{conregint} we conclude that $\tilde{X}_H$ is singular.

Hence, to prove the existence of a
singular fiber we may assume that $X$ is nonsingular.
After a linear transformation, we may assume that $L\subset \P^3$ is
given by $w=z=0$ and $X$ is given by $F=0$ for some homogeneous
polynomial $F\in k[x,y,z,w]$ of degree $3$. Let $P\in X 
\subset \P^3$ be a point where both $\partial F/\partial x$ and
$\partial F/\partial y$ vanish (the existence of $P$ follows from the
Projective Dimension Theorem, see \myref{\hag}, Thm. I.7.2). 
Set $t_0=(\partial F/\partial z)(P)$ and
$t_1=(\partial F/\partial w)(P)$. Then $t_0$ and $t_1$ are not both zero
because $P$ is nonsingular. The tangent plane $T_{P}$ to $X$ at $P$ is then
given by $t_0z+t_1w=0$, so it contains $L$. The fiber $T_P \cap X$
above $T_P$ is singular, as $T_P$ is tangent at $P$.
\end{proof}

\begin{lemma}\label{ratsurffour}
Under the assumptions of Proposition \ref{rationalellipticsurface}
the morphism $\tilde{f}$ is a relatively minimal fibration.
\end{lemma}
\begin{proof}
By Lemmas \ref{ratsurfone}, \ref{ratsurfonehalf}, and
\ref{ratsurftwo}, the hypotheses of 
Theorem \ref{canonicalonminimalfibration} are satisfied, so it
suffices to show that $K_{\tilde{X}}^2=0$. 
Let $\rho\colon X'\rightarrow X$ be the blow-up of $X$ at the 
three points $M_i$, and let $\sigma \colon \tilde{X} \rightarrow X'$
be the minimal desingularization of $X'$. 

For any projective variety
$Z$, let $K_Z^\circ$ denote the divisor associated to the dualizing
sheaf $\omega_Z^\circ$, see \myref{\hag}, \parag{} III.7. If $Z$ is
nonsingular, then $K_Z^\circ$ is linearly equivalent to the
canonical divisor $K_Z$, see \myref{\hag}, Cor. III.7.12.
From \myref{\hag}, Thm. III.7.11, we find that 
$\omega_X^\circ \isom \O_X(d-4)$ with $d = \deg X = 3$. Hence, if $H$
is a hyperplane that does not meet any of the $M_i$ or the singular
points of $X$, then $K_X^\circ$ is linearly equivalent to $-(H \cap X)$.

Let $U$ be the maximal smooth open subset of $X$, and set 
$V=\rho^{-1}(U)$. As $V$ is isomorphic to
$U$, blown up at three nonsingular points, we find by \myref{\hag},
Prop. V.3.3, that $K_V = \rho^* K_U + \M_1+\M_2+\M_3$. Since $\rho$
is an isomorphism outside the $M_i$, we find that $K_{X'}^\circ =
\rho^*K_X^\circ +\M_1+\M_2+\M_3$. As
$K_X^\circ$ does not meet the $M_i$, and $\M_i^2=-1$ (see
\myref{\hag}, Prop. V.3.2) we get 
$$
(K_{X'}^\circ)^2 = 
(\rho^*K_X^\circ)^2+\M_1^2+\M_2^2+\M_3^2 = (K_X^\circ)^2-3 = (H\cap
X)^2-3 = \deg X -3 =0.
$$
Du Val \myref{\duval} proves that rational singularities do not
affect adjunction, i.e., there is an isomorphism 
$\omega_{\tilde{X}}^\circ \isom 
\sigma^* \omega_{X'}^\circ$, see also \myref{\pink}, \parag{} 15,
Prop. 2, and \parag{} 17. \donecheck{ask John for reference of Artin}
Hence, we get $K_{\tilde{X}} \sim
K_{\tilde{X}}^\circ \sim \sigma^* K_{X'}^\circ$. As $\sigma$ is an
isomorphism on $\sigma^{-1}(V)$, we get  
$K_{\tilde{X}}^2 = (\sigma^* K_{X'}^\circ)^2 = (K_{X'}^\circ)^2=0$.
%
%
\end{proof}

\begin{proofof}{\bf Proposition \ref{rationalellipticsurface}.}
This follows immediately from Lemmas \ref{ratsurfone},
\ref{ratsurfonehalf}, \ref{ratsurftwo}, \ref{ratsurfthree}, 
and \ref{ratsurffour}.
\end{proofof}

\begin{proofof}{\bf Proposition \ref{ratsurfez}.} Follows immediately
  from Proposition \ref{rationalellipticsurface}.
\end{proofof}

\begin{remark}
If $L$ intersects $X$ in one of its singular points, then one
could still define a fibration $\tilde{X} \rightarrow \P^1$ in the
same way as in Proposition \ref{rationalellipticsurface}. 
For almost all hyperplanes $H$ the fiber above $H$ will be the
normalization of the singular cubic curve $H \cap X$.
Hence this will not be an elliptic fibration.
\end{remark}

\begin{remark}
In characteristic $3$, all fibers might be singular, as is the case when
$X$ is given by $y^2z+yz^2+wxy+wxz+xz^2+wy^2=0$ and $L$ is given by
$x=w=0$. The intersection of $X$ with the plane $H_t$ given by $w=tx$ is
singular at the point $[x\mcol y \mcol z \mcol w] = [1\mcol
  t^{1/3} \mcol t^{2/3} \mcol t]$ on the twisted cubic curve in $\P^3$. 
The plane $H_t$ is tangent to $X$ at that point. The only singular
points of $X$ are three ordinary double points at $[1:0:0:0]$,
$[0:0:0:1]$, and $[1:1:1:1]$.

In characteristic $2$, we can also get all fibers to be singular, as
one easily checks in case $X$ is given by $x^3+x^2z+x^2w+y^3+yzw=0$
and $L$ is given by $w=z=0$. The only singular points on $X$ are 
the ordinary double points $[0:0:0:1]$ and $[0:0:1:0]$. 

In the proof of Proposition \ref{rationalellipticsurface} the fact
that the characteristic of $k$ is not equal to $2$ or $3$ is only used
in Lemma \ref{ratsurftwo}. Hence the conlusion
of the proposition is also true in characteristic $2$ and $3$ if we
add to the hypotheses that almost all planes through $L$ are not
tangent to $X$. By Bertini's Theorem, the set of planes that intersect
$X$ in a nonsingular curve is open (see \myref{\hag}, Thm. II.8.18), 
so it is enough to require that there is at least one plane through
$L$ that is not tangent to $X$.
\end{remark}

\begin{remark}
The singular points on $X$ as in Proposition \ref{rationalellipticsurface}
can be used to find sections of $\tilde{f}$. If $X$ has two
singular points $P$ and $Q$, then the line $l$ through $P$ and $Q$ lies on
$X$, for if it did not, it would have intersection multiplicity at
least $4$ with $X$, but by B\'ezout's Theorem the intersection
multiplicity should be $3$. Therefore, either $l$ intersects $L$ and 
thus $l$ is contained in the fiber
above the plane that $L$, $P$, and $Q$ all lie in, or $l$ gives a
section of $\tilde{f}$.
\end{remark}

The next proposition describes how to construct an elliptic surface by
base extension of another elliptic surface. This construction will
also be used in the proof of Theorem \ref{maintheorem}.

\begin{proposition}\label{baseextendellipticsurface}
Let $S$ be an elliptic surface over a smooth, irreducible, projective
curve $C$ over an algebraically closed field $k$, with fibration 
$g$ and section $\O$ of $g$.  Let $\gamma\colon C'\rightarrow C$ 
be a nonconstant map of curves from a smooth,
irreducible, projective curve $C'$, which is unramified above those
points in $C$ where $g$ has singular fibers. 
Set $S' = S\times_C C'$, let $g'$ be the projection $S' \ra C'$,  
and let $\O'\colon C' \ra S'$ denote the morphism induced by the
identity on $C'$ and the composition $\O \circ \gamma$. Then 
$\O'$ is a section of $g'$ and they make $S'$ into an elliptic surface
over $C'$. The Euler characteristics $\chi_S=\chi(\O_S)$ and 
$\chi_{S'}=\chi(\O_{S'})$ are related by $\chi_{S'} = (\deg \gamma)
\chi_S$. 
\end{proposition}

$$
\xymatrix{
S' \ar @<+.2mm>[dd]^{g'}\ar[rr]&& S \ar @<+.2mm>[dd]^{g} \\
&& \\
C'\ar @/^/ @<+.4mm> [uu]^{\O'} \ar[rr]_{\gamma}&& C\ar @/^/ @<+.4mm> [uu]^{\O}
}
$$

\begin{proof}
Since projective morphisms are stable under base extension (see
\myref{\hag}, exc. II.4.9), we find that $S'$ is projective over $C'$,
which is projective over $\Spec k$, so $S'$ is projective. The
composition $g' \circ \O'$ is by construction the identity on $C'$, so
$\O'$ is a section of $g'$.

As $k$ is algebraically closed, the residue field $k(x)$ of a closed
point $x \in C'$ is isomorphic to the residue field $k(\gamma(x))$. 
Hence the fiber above $x$ is isomorphic
to the fiber above $\gamma(x)$, as we have 
$$
\Spec k(x) \times_{C'} S' \isom \Spec k(x) \times_{C'} C' \times_C S
\isom \Spec k(x) \times_C S \isom \Spec k(\gamma(x)) \times_C S.
$$
Therefore, as for $g$, all fibers of $g'$ are connected. As $g$ is elliptic, 
all but finitely many fibers of $g'$ will be smooth curves of genus
$1$. Since $g$ has a singular fiber, so does $g'$. 
From Lemma \ref{connfibers} we find that $g'_* \O_{S'} \isom \O_{C'}$.
As $C'$ is irreducible and projective, this implies
$\dim H^0(S', \O_{S'}) = \dim H^0(C',g'_* \O_{S'})=
\dim H^0(C',\O_C') = 1$. We conclude that $S'$ is connected.

To prove that $S'$ is smooth and irreducible, set 
$h = \gamma \circ g'$. By assumption
there are open $U, V \subset C$
with $U \cup V = C$, such that $\gamma|_{\gamma^{-1}(U)}$ is
unramified, whence smooth, and $g|_{g^{-1}(V)}$ has no singular
fibers, which implies it is smooth by Remark \ref{smooth}. As smooth
morphisms are stable under base extension and composition (see
\myref{\hag}, Prop. II.10.1), we find first that $h^{-1}(U) = 
g^{-1}(U) \times_U \gamma^{-1}(U)$ is smooth over $g^{-1}(U)\subset
S$. As $S$ is smooth over $k$ and $g^{-1}(U)$ is open in $S$, 
we conclude that $h^{-1}(U)$ is smooth
over $k$. Similarly, $h^{-1}(V)$ is smooth over $k$, whence so is
their union $S'$.
%
%
%
As $S'$ is also connected, we find that $S'$ is irreducible from Lemma
\ref{conregint}. 

To prove that $g'$ is relatively minimal, it suffices 
by Theorem \ref{blowdownable} to show that no fiber $S'_x$ above $x\in C'$
contains an exceptional prime divisor. 
Let $D'$ be an irreducible component of the fiber $S'_x$, mapping
isomorphically to the irreducible component $D$ of $S_{\gamma(x)}\isom
S'_x$ under the induced morphism $\gamma' \colon S' \rightarrow S$. 
Suppose that $D'$ is an exceptional divisor, i.e., $D' \isom \P^1$ and
$D'^2=-1$. If
$\gamma(x)$ is contained in $V$, then the fiber $S_{\gamma(x)}$, 
and hence $S'_x$, is smooth. As all fibers are connected, $S'_x$ is then
irreducible, so $D'=S'_x$. Since any fiber is numerically equivalent to any
other, this implies $D'^2=0$, contradiction.
Therefore, we may assume that $\gamma(x) \not \in V$, so $\gamma(x)
\in U$ and $D' \subset h^{-1}(U)$. As \'etale morphisms are stable under
base extension and $\gamma|_{\gamma^{-1}(U)}$ is \'etale, we find that
$\gamma'|_{h^{-1}(U)}$ is \'etale.

%
 
For any morphism of schemes $\varphi \colon X \rightarrow Y$, let
$\Omega_{X/Y}$ denote the sheaf of relative differentials of $X$ over
$Y$. If $X$ is a nonsingular variety over $k$, then let $\T_X$ denote
the tangent sheaf $\shfHom(\Omega_{X/k}, \O_X)$. For any nonsingular
subvariety $Z \subset X$, let $\Normal_{Z/X}$ denote the normal sheaf
of $Z$ in $X$, see \myref{\hag}, \parag{} II.8.

We will show that the self-intersection number $D'^2 = \deg
\Normal_{D'/S'}$ on $S'$ (see \myref{\hag}, example V.1.4.1) 
is equal to the self-intersection number $D^2 =
\deg \Normal_{D/S}$. Since $D$ is not an exceptional curve, that implies
that $D'^2 \neq -1$, which is a contradiction.
As $\gamma'$ induces an isomorphism from $D'$ to $D$, it suffices to
show that $\Normal_{D'/S'}$ is isomorphic to $\gamma'^*
\Normal_{D/S}$.

There is an exact sequence 
\begin{equation}\label{sesone}
0 \ra \T_{D'} \ra \T_{S'} \otimes \O_{D'} \ra \Normal_{D'/S'} \ra 0
\end{equation}
(see \myref{\hag}, page 182),
and by applying the isomorphism $(\gamma'|_{D'})^*$ to the similar
sequence for $D$ in $S$ we also get the exact sequence
\begin{equation}\label{sestwo}
0 \ra\gamma'^* \T_D \ra\gamma'^*(\T_S \otimes \O_{D}) \ra
\gamma'^* \Normal_{D/S} \ra 0.
\end{equation}
The natural morphisms $\T_{D'} \ra \gamma'^* \T_D$ and 
$\T_{S'} \otimes \O_{D'}\ra \gamma'^*(\T_S \otimes \O_{D})$ induce a
morphism between the short exact sequences (\ref{sesone}) and 
(\ref{sestwo}). To prove that the last morphism 
$\Normal_{D'/S'} \ra\gamma'^* \Normal_{D/S}$
is an isomorphism, it suffices by the snake lemma 
to prove that the first two are. 
Clearly,  $\T_{D'} \ra \gamma'^* \T_D$ is an isomorphism of sheaves on
$D'$, as $\gamma'|_{D'}$ is an isomorphism. To show that 
$$
\T_{S'} \otimes \O_{D'}\ra \gamma'^*(\T_S \otimes \O_{D}) \isom
   \gamma'^* \T_S \otimes \gamma'^* \O_D 
                      \isom \gamma'^* \T_S \otimes \O_{D'}
$$ 
is an isomorphism, it suffices to show that $\T_{S'} \ra \gamma'^*
\T_S$ is an isomorphism on the open subset $h^{-1}(U)\subset S'$ 
containing $D'$. 
This is true, as by \myref{\sgaone}, Expos\'e II, Cor. 4.6,
a morphism $f \colon X \ra Y$ of smooth
$T$-schemes is \'etale if and only if the morphism 
$f^* \Omega_{Y/T} \ra \Omega_{X/T}$ is an isomorphism. Taking the dual
gives what we need, if we choose $T=\Spec k$, and
$f=\gamma'|_{h^{-1}(U)}$. 

For the last statement we will use that 
by \myref{\kodtwo}, Thm. 12.2, we have
\begin{equation}\label{eulerchiS}
\eqalign{
12\chi_S=\mu+6\sum_{b\geq 0}\nu(I_b^*)&+2\nu(II)+10\nu(II^*)+3\nu(III) \cr 
&+9\nu(III^*)+4\nu(IV)+8\nu(IV^*),
}
\end{equation}
where $\nu(T)$ is the number of singular fibers of $g$ of type $T$ and $\mu$
is the degree of the map $j(S/C) \colon C \rightarrow \P^1$, sending
every element $x\in C$ to the $j$-invariant of the fiber $S_x$.

As the morphism $\gamma$ is unramified above
the points of $C$ where $g$ has singular fibers, it follows
that the singular fibers of $g'$ come in $n$-tuples, with $n=\deg
\gamma$. Each $n$-tuple consists of
$n$ copies of one of the singular fibers of $g$. 
Hence, if $\nu'(T)$ denotes the number of singular fibers of $g'$ of
type $T$, then we have $\nu'(T) = n\nu(T)$. As $j(S'/C')$ is the
composition of $\gamma$ and $j(S/C)$, we also get $\mu' = n\mu$, where
$\mu'$ is the degree of $j(S'/C')$. From (\ref{eulerchiS}) and its
analogue for $S'$ we conclude that $\chi_{S'} = n \chi_S$.
\end{proof}

\bigskip

\mysection{\titlefive}\label{secfive}

Let $X$, $G$, and $f \colon X \dashrightarrow \P^1$ be as in section 2.
The rational map $f$ is defined everywhere, except at the three
intersection points $M_1=[0:0:1:-1]$, $M_2=[0:1:0:-1]$, and $M_3=[0:1:-1:0]$
of $X$ with the line $L$ given by $r=x+y+z=0$. 
A straightforward computation shows that $X$ has three singular points 
$N_1=[0:1:0:0]$, $N_2=[0:0:1:0]$, and $N_3=[0:0:0:1]$, all ordinary
double points, forming a full orbit 
under $G$, and all mapping to $[0:1]$ under $f$.  
Let $\pi \colon \tilde{X} \ra X$ be the blow-up of $X$ at the six points $M_i$ 
and $N_i$. Let $\M_i$ and $\N_i$ denote the exceptional curves above $M_i$ and 
$N_i$ respectively.

\begin{proposition}\label{Xelliptic}
The surface $\tilde{X}$ is smooth. 
The rational map $f \circ \pi$ extends to a
morphism $\tilde{f} \colon \tilde{X} \rightarrow \P^1$.  
It maps the $\tilde{M}_i$ isomorphically to $\P^1$ and together with
the section  
$\O = \tilde{f}|_{\tilde{M}_3}^{-1}$ it makes $\tilde{X}_k$ into an
elliptic surface over $\P^1$ for any algebraically closed field $k$ of
characteristic $0$.
\end{proposition}
\begin{proof}
Ordinary double points are resolved by blowing up once, see
\myref{\hag}, exc. I.5.7. Hence $\tilde{X}$ is the minimal
desingularization of $X$ blown up at the $M_i$.
The rational map $f$ sends all points of $X$ 
(except for the $M_i$) in the plane through $L$ given by 
$t_1r=t_0(x+y+z)$ to the point $[t_0:t_1]$. 
Hence this proposition follows from
Proposition \ref{rationalellipticsurface}.
\end{proof}

\begin{remark}\label{singfib}
In this explicit case, it would have been easier to check by hand that
$\tilde{f}$ makes $\tilde{X}_k$ into an elliptic surface over
$\P^1$. From Theorem \ref{blowdownable} it follows that, in order to
prove that $\tilde{f}$ is a minimal fibration, it suffices to check that no
reducible fiber contains a rational curve with self-intersection $-1$.
As the only singular points of $X$ lie above $[0\mcol 1]\in \P^1$,
it follows that for all $\tnum \neq 0, \infty$, the fiber 
$\tilde{X}_{\tnum}$ above $[\tnum \mcol 1]$ is
given by the intersection of $X$ with the plane given by
$r=\tnum(x+y+z)$. Hence for $\tnum \neq 0,\infty$, the fiber 
is isomorphic to the plane curve given by $\tnum^2(x+y+z)^3=xyz$, which is 
nonsingular as
long as $\tnum(27\tnum^2-1) \neq 0$. For $\tnum$ with $27\tnum^2=1$ we get a
nodal curve, whence 
a fiber of type $I_1$, following the Kodaira-N\'eron classification of
special fibers, see \myref{\silvtwo}, IV.8 and \myref{\kod}. At $\tnum=0$ and
$\tnum=\infty$ one checks that the fibers are of type $I_6$ and $IV$
respectively. None of these fibers contains an exceptional curve.
\check{$\tnum \rightarrow t$?}
\end{remark}

\begin{remark}\label{noparamsfixedt}
From the previous remark, it follows that the fiber of $\tilde{f}$
above every {\em rational} point $[\tnum\mcol 1] \in \P^1$ with
$\tnum>0$, is a 
curve of genus $1$, which can therefore not be rationally
parametrized. Therefore, there is no rational parametrization of infinitely
many rational triangles, all having the same area and the same perimeter.
\check{$\tnum \rightarrow t$?}
\end{remark}

\begin{remark}
Later we will see a Weierstrass form for the generic fiber of
$\tilde{f}$. Based on that, Tate's algorithm (see
\myref{\silvtwo}, IV.9 and \myref{\tate}) describes the
special fibers of a minimal proper regular model. They coincide with
the fibers 
described in Remark \ref{singfib}, which gives another proof of the fact that
$\tilde{f}$ is relatively minimal.
\end{remark}

Let $E$ denote the generic fiber of $\tilde{f}$, an elliptic curve
over $k(\P^1)\isom \Q(t)$. By Lemma \ref{sections} we can identify the 
sets $\tilde{X}(\P^1)$ and $E(k(\P^1))$.
The curve $E$ is isomorphic to the plane curve in $\P^2_{\Q(t)}$ given by 
\begin{equation}\label{firstmodel}
t^2(x+y+z)^3 = xyz.
\end{equation}
The origin $\O = \M_3$ then has coordinates $[x\mcol y\mcol z] = 
[1 \mcol -1 \mcol 0]$. Let $P$ denote the section 
$\M_1 = [0\mcol  1\mcol -1]$. 
A standard computation shows that the $\M_i$ correspond with
inflection points. As they all lie on the line given by $x+y+z=0$, 
we find that $P$ has order $3$ and $2P=\M_2=[1\mcol 0 \mcol -1]$.
This also follows from the following lemma, which 
gives a different interpretation of the action of $G$.

\begin{lemma}\label{automs}
The automorphism $\tilde{X} \rightarrow \tilde{X}$
induced by the $3$-cycle $(x\,y\,z)$ on the coordinates of $X$ corresponds
with translation by $P$ on the nonsingular fibers and the generic
fiber of $\tilde{f}$. Similarly, the automorphism induced by $(x\,y)$
corresponds with multiplication by $-1$. 
\end{lemma}
\begin{proof}
The automorphism $\psi$ induced by $(x\,y)$ fixes $\O$, so it induces an
isogeny on the nonsingular fibers. The automorphism group of an
elliptic curve over $\C$ is isomorphic to the group of roots of unity in 
$\Z$ or in a quadratic
order. As we have $\psi^2=1$ and $\psi \neq 1$, we find $\psi =
[-1]$. 

Let $\varphi = T_{-P} \circ (x\,y\,z)$ be the composition of translation
by $-P$ and the automorphism induced by $(x\,y\,z)$. Then $\varphi$ fixes
$\O$, so it induces an isogeny as well. The cube of $(x\,y\,z)$ then sends
any point $Q$ to $P+\varphi(P)+\varphi^2(P)+\varphi^3(Q)$. As
$(x\,y\,z)^3=1$, we find that $\varphi^3=1$. The involution
$(y\,z)=(x\,y)(x\,y\,z)$ sends $Q$ to $-P-\varphi(Q)$, so $(y\,z)^2$ sends $Q$
to $-P+\varphi(P)+\varphi^2(Q)$. Since $(y\,z)^2=1$, we also find that
$\varphi^2=1$, so $\varphi=1$. Hence $(x\,y\,z) = T_P$ on $E$. As $E$ is
dense in $\tilde{X}$, we find $(x\,y\,z) = T_P$ on $\tilde{X}$, see 
\myref{\hag}, exc. II.4.2.  
\end{proof}

As mentioned before, we want infinitely many $\tnum$ for which the fiber
$X_{\tnum}$ above $[\tnum\mcol 1]$ has infinitely many rational points
$[r_i\mcol x_i \mcol y_i \mcol z_i]$ with $x_i/r_i$, $y_i/r_i$,
$z_i/r_i > 0$, and all in different orbits under $G$. 
If the Mordell-Weil rank of $E(\Q(t))\isom \tilde{X}(\P^1)$ had been 
positive, we might
have been able to find infinitely many such points for almost all 
rational $\tnum$ satisfying some inequalities. 
Unfortunately, the next theorem tells us that this is not the case.
\check{$\tnum \rightarrow t$? think no}

\begin{theorem}\label{mordellrat}
The Mordell-Weil group $E(\C(t))$ is isomorphic to $\Z \times \Z/3\Z$.
It is generated by the $3$-torsion point $P$ and the point $Q\colon 
[r:x:y:z] = [t:it:-it:1]$. The Mordell-Weil group $E(\R(t))$ is 
equal to $\langle P\rangle \isom \Z/3\Z$.
\end{theorem}
\begin{proof}
As $\tilde{X}$ is rational, the N\'eron-Severi group $\NS(\tilde{X}_\C)$ 
is a unimodular lattice of rank $10$, see \myref{\shioda}, Lemma 10.1. 
Let $T \subset \NS(\tilde{X}_\C)$ be as in Proposition \ref{Tkern}. From
Remark \ref{singfib} and Proposition \ref{Tkern}, we find that $T$ has rank
$2+(6-1)+(3-1)+(1-1)+(1-1)=9$ and we can find explicit generators. 
Consider the lattice $T+\langle (P),(Q)\rangle$. 
Computing the explicit intersections of our
generators, we find that the lattice $T+\langle (P),(Q)\rangle$ has rank 10,
whence finite index in $\NS(\tilde{X}_\C)$. Also, it is already
unimodular, whence equal to $\NS(\tilde{X}_\C)$. Hence, $E(\C(t))$ is
generated by $P$ and $Q$ and has rank $1$.

Complex conjugation on $Q$ permutes the $x$- and $y$-coordinates, so by Lemma
\ref{automs} we find $\overline{Q}=-Q$ in $E(\C(t))$. 
If $mQ+nP$ is real for some
integers $m,n$, then so is $mQ$ and hence $mQ=m\overline{Q}=-mQ$, so $2mQ=0$. 
Since $Q$ has infinite order, we conclude that $m=0$, so $E(\R(t)) = \langle P
\rangle$.
\end{proof}

%
%
%

To find more curves over $\Q$, we will apply a base change to our base curve
$\P^1$ by a rational curve on $\tilde{X}$. As we have a
parametrization of $X$, it is
easy to find such a curve. Taking $u=\svar$ and $v=w=1$ in
(\ref{Xrational}) we find a curve $C$ on $X$ parametrized by
$$
\beta \colon \P^1 \rightarrow C \colon [\svar \mcol 1] \mapsto
[r\mcol x\mcol y \mcol z ]=[\svar-1 \mcol \svar+1 \mcol \svar-1
  \mcol \svar(\svar-1)].
$$
We will denote its strict transform on $\tilde{X}$ by $C$ as well.
The map $\tilde{f}$ induces a $2$-$1$ map from $C$ to $\P^1$. 
The composition $\tilde{f} \circ \beta$ is given by $[\svar\mcol 1]
\mapsto [\svar-1 \mcol \svar(\svar+1)]$. Hence, if we identify the
function field $K=k(C)$ of $C$ with $\Q(\svar)$, then the field
extension $K/k(\P^1)$ is given by $\Q(t) \hookrightarrow \Q(\svar)
\colon t \mapsto (\svar-1)\svar^{-1}(\svar+1)^{-1}$.
Thoughout the rest of this article, as in Theorem \ref{maintheorem}
and Remarks \ref{singfib} and \ref{noparamsfixedt}, one should think
of $\snum$ and $\tnum$ as specific values for the indeterminates
$\svar$ and $t$ respectively, unless the context clearly suggests
otherwise. \check{say earlier?} 

Let $Y$ denote the fibered product  $\tilde{X} \times_{\P^1} C$, let
$\delta$ denote the projection $Y \rightarrow \tilde{X}$, and
let $g$ denote the projection $Y \rightarrow C$. The generic fiber of
$g$ is isomorphic to $E_K = E \times_{k(\P^1)} K$.  
The identity on $C$ and the composition $\O \circ\tilde{f}|_C \colon C
\ra \tilde{X}$ together induce a section $C \ra Y$ of $g$, which we
will also denote by $\O$. The closed immersion $C \ra \tilde{X}$ and
the identity on $C$ together induce a section of $g$ which we will
denote by $R$. 

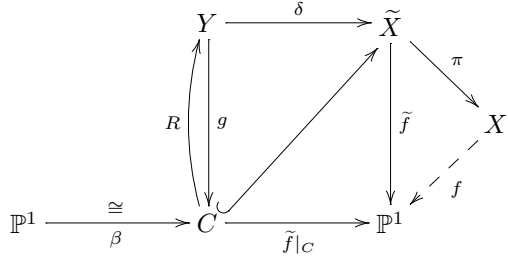
$$
\xymatrix{&& Y\ar @<+.2mm>[dd]^g\ar[rr]^{\delta}&&\tilde{X}\ar[dr]^{\pi}
                                            \ar[dd]^{\tilde{f}}& \\
          && & & & X \ar @{-->}[ld]^{f} \\
         \P^1 \ar[rr]^{\isom}_\beta && C \ar @/^/ @<+.4mm> [uu]^R
            \ar[rr]_{\tilde{f}|_C}\ar @<-1mm> @{^{(}->}[uurr]&& \P^1&
}
$$

\begin{proposition}
The fibration $g$ and its section $\O$ make $Y_k$ into an elliptic
surface over $C_k$ for any algebraically closed field $k$ of
characteristic $0$.
\end{proposition}
\begin{proof}
One easily checks that $\tilde{f}|_C \colon C \rightarrow \P^1$ is
unramified at the points of $\P^1$ where $\tilde{f}$ has singular fibers.
Hence, this proposition follows immediately from Proposition
\ref{baseextendellipticsurface} and Proposition \ref{Xelliptic}. 
\end{proof}

From (\ref{firstmodel}) we find that $E_K$ is isomorphic to the plane
cubic over $K$ given by 
$$
(\svar-1)^2(x+y+z)^3 = \svar^2(\svar+1)^2xyz.
$$
The linear transformation 
\begin{equation}\label{transform}
p=-4(\svar-1)^2(x+y)z^{-1}, \qquad q = 4(\svar-1)^2\svar(\svar+1)(x-y)z^{-1},
\end{equation}
or, equivalently,
\begin{equation}\label{trafoback}
\eqalign{
x &= -\svar(\svar+1)p+q, \cr
y &= -\svar(\svar+1)p-q, \cr
z &= 8(\svar-1)^2 \svar (\svar+1),
}
\end{equation}
gives the Weierstrass equation
\begin{equation}\label{weier}
q^2=(p-4(\svar-1)^2)^3+\svar^2(\svar+1)^2p^2=F(\svar,p).
\end{equation}
with
\begin{equation}\label{jdisc}
\eqalign{
j&=j(E_K)=j(E)=\frac{(24t^2-1)^3}{t^6(27t^2-1)}, \cr
\Delta&=2^{12}(\svar-1)^6\svar^4(\svar+1)^4(\svar^4+2\svar^3-26\svar^2+54\svar-27).
}
\end{equation}

The Weierstrass coordinates of $P$ and $R$ are given by
$$
\eqalign{
(p_P,q_P)&=(4(\svar-1)^2,4\svar(\svar+1)(\svar-1)^2) \qquad \text{and} \cr
(p_R,q_R)&=(8-8\svar,8\svar^2-8).\cr
}
$$

\begin{lemma}\label{finish}
The section $R$ has infinite order in the group $Y(C) \isom E_K(K)$. 
\end{lemma}
\begin{proof}
 The $p$-coordinate of $2R+P$ equals 
$4(\svar^4-6\svar^3+10\svar^2-2\svar+1)(\svar-1)^{-2}$, so $2R+P$ 
is contained in the
kernel of reduction at $\svar-1$. In characteristic $0$ the kernel of
reduction has no nontrivial torsion (see \myref{\silv},
Prop. VII.3.1), so we find that $2R+P$ has infinite order, whence so
does $R$.
\end{proof}

For every
integer $n \geq 1$, let $\gamma_n \colon \P^1 \ra X$ denote the composition 
\begin{equation}\label{composition}
\P^1 \stackrel{\beta}{\longrightarrow} C
\stackrel{(2n-1)R}{\longrightarrow}
Y \stackrel{\delta}{\longrightarrow} \tilde{X}
\stackrel{\pi}{\longrightarrow} X. 
\end{equation}

Theorem \ref{maintheorem} will follow from the following
proposition.

\begin{proposition}
Let $\snum>1$ be a rational number. For every integer $n \geq 1$, let
$r_n,x_n,y_n$, and $z_n$ be such that
$\gamma([\snum:1])=[r_n:x_n:y_n:z_n]$ and set 
$$
a_n=\frac{(\snum-1)(y_n+z_n)}{r_n}, \quad 
b_n=\frac{(\snum-1)(x_n+z_n)}{r_n}, \quad 
c_n=\frac{(\snum-1)(x_n+y_n)}{r_n}.
$$
Then for every $n\geq 1$ there is a triangle $\Delta_n$ with sides
$a_n$, $b_n$, $c_n$, perimeter $2\snum(\snum+1)$, inradius 
$\snum-1$, and area $\snum(\snum^2-1)$. The triangles $\Delta_n$
are pairwise nonsimilar.
\end{proposition}
\begin{proof}
Let a rational $\snum>1$ be given and set $c =
\beta([\snum:1])\in C$. Then $\tilde{f}|_C(c) = [\tnum:1]$ for $\tnum =
(\snum-1)\snum^{-1}(\snum+1)^{-1}>0$, so the fiber $Y_c$ is
isomorphic to the fiber $\tilde{X}_{\tnum}$ of $\tilde{f}$ above $[\tnum:1]$.
By Remark \ref{singfib}, this
fiber is nonsingular and isomorphic to the intersection $E_{\tnum}$ of
$X$ with the hyperplane given by $r=\tnum(x+y+z)$. This intersection
$E_{\tnum}$ can be given the structure of an elliptic curve with $M_3$
as origin. The specialization map $Y(C) \rightarrow Y_c(\Q) \colon 
S \mapsto S \cap Y_c=S(c)$ induces a homomorphism
$\psi \colon Y(C) \rightarrow E_{\tnum} \subset X$ sending a section $S$
of $g$ to $\pi(\delta(S(c)))$. Set $\Theta_n = \gamma_n([\snum:1])\in
X= [r_n:x_n:y_n:z_n]$.
Then we have $\Theta_n = \psi((2n-1)R)\in E_{\tnum}$, so on $E_{\tnum}$ we get
$\Theta_n = (2n-1) \Theta_1$. The elliptic curve $E_{\tnum}$ has a
Weierstrass model $q^2=F(\snum,p)$, see (\ref{weier}). For $n \geq 1$, let
$(p_n,q_n)$ denote the Weierstrass coordinates of $\Theta_n$, so 
$(p_1,q_1) = (8-8\snum,8\snum^2-8)$. 
\check{$\tnum \rightarrow t$? think no}

Note that $F(\snum,0)=-64(\snum-1)^6<0$, but for $p_1=8-8\snum<0$ we have 
$F(\snum,p_1)=q_1^2>0$. 
We conclude that for any real point on $E_{\tnum}$ with Weierstrass
coordinates $(p,q)$, the condition $p<0$ is equivalent to the point lying on
the real connected component of $E_{\tnum}$ that does not contain $\O$.
Since $\Theta_1$ lies on this component, so do all its odd multiples
$\Theta_n$.

To prove that there is a triangle with sides $a_n$, $b_n$, and $c_n$,
it suffices to prove $x_n/r_n, y_n/r_n, z_n/r_n >0$. 
If $\Theta_n=M_i$ for $i=1,2$, or $3$, then
$3\Theta_n=\O$, which contradicts the fact that $\Theta_n$ lies on the
real component of $E_{\tnum}$ that does not contain $\O$. Hence
$f$ is well-defined at $\Theta_n$ and from
$[r_n:x_n+y_n+z_n]=f(\Theta_n) = [\tnum:1]$, with $\tnum>0$, 
we find $r_n\neq 0$ and
$x_n+y_n+z_n \neq 0$, whence $x_ny_nz_n \neq 0$. To make computations
easier, we may assume that $z_n = 8(\snum-1)^2\snum(\snum+1)>0$. 
We will show that then also $x_n,y_n,r_n >0$. By (\ref{trafoback}) we get 
$$
\eqalign{
x_n &= -\snum(\snum+1)p_n+q_n, \cr
y_n &= -\snum(\snum+1)p_n-q_n. \cr
}
$$
Hence the condition $x_n,y_n>0$ is equivalent to 
$$
p_n<0 \qquad \mbox{and} \qquad q_n^2<\left(\snum(\snum+1)p_n\right)^2=
q_n^2-(p_n-4(\snum-1)^2)^3.
$$
The last inequality is equivalent to
$p_n<4(\snum-1)^2$, which is automatically satisfied if the first
inequality $p_n<0$ is satisfied. 
Therefore, the condition is equivalent to the inequality
$p_n<0$, which is satisfied as it is equivalent to $\Theta_n$ lying on
the real connected component that does not contain $\O$.
From $r_n = \tnum(x_n+y_n+z_n)$ we also conclude $r_n>0$. 
The triangle with sides $a_n$, $b_n$, and $c_n$ has inradius 
$\snum-1$, perimeter
$2(\snum-1)(x_n+y_n+z_n)/r_n = 2(\snum-1)/\tnum = 2\snum(\snum+1)$
and hence area $\snum(\snum^2-1)$. 

In order to prove that all the $\Theta_n$ are different, assume that
$\Theta_1$ has finite order. As $\Theta_1$ lies on the real component
that does not contain $\O$, it has even order, so by Mazur's Theorem
(see \myref{\silv}, Thm. III.7.5 for statement, 
\myref{\maz}, Thm. 8 for
the proof) we find that $m\Theta_1=\O$ for $m=8,10$, or $12$.
For each of these three values for $m$ we can compute explicit 
rational functions $\xi_m,\eta_m \in \Q(\svar)$ such that the 
coordinates of $m\Theta_1$ are given by $(\xi_m(\snum),\eta_m(\snum))$.
For $m=8,10$, or $12$, these rational functions turn out to not have
any rational poles, so $\Theta_1$ has infinite order.

To show that the triangles are pairwise nonsimilar, it suffices by
Lemma \ref{TrianglesToPoints} to show that the $\Theta_n$ lie in
different orbits under 
$G$. Suppose that $\Theta_n$ and $\Theta_{n'}$ are in the same orbit
under $G$ for some $n,n' \geq 1$. Then by Lemma \ref{automs} we get
$\Theta_n = \pm \Theta_{n'}+kP$ for $k=0,1$ or $2$.  
Hence $3\left((2n-1)\mp (2n'-1)\right)\Theta_1=
3(\Theta_n\mp \Theta_{n'})=3kP=\O$, 
so $2n-1 = \pm (2n'-1)$, as $\Theta_1$ has
infinite order. From $n,n' \geq 1$ we find $n=n'$ and hence $k=0$.
Thus, $\Theta_n = \Theta_{n'}$.
\end{proof}

\begin{proofof}{\bf Theorem \ref{maintheorem}.}
Consider the open affine subset $U\subset X$ defined by $r \neq 0$, 
which is isomorphic to $\Spec A$ for 
$A=\Q[x,y,z]/\big(x+y+z-xyz\big)$. For each $n \geq 1$,
let $V_n \subset \P^1$ be a dense open affine subset such that 
the composition $\gamma_n$ of morphisms in (\ref{composition}) maps $V_n$ to
$U$. This is possible because the image of $\P^1$ is not entirely
contained in the closed subset of $X$ given by $r=0$.
Then there is a ring $B_n \subset \Q(\svar)$ such that $V_n$ is
isomorphic to $\Spec B_n$ and the composition in (\ref{composition}) is
given by a ring homomorphism 
$\varphi_n \colon A \rightarrow B_n \subset \Q(\svar)$. 
Let $x_n(\svar)$, $y_n(\svar)$, $z_n(\svar) \in \Q(\svar)$ 
be the images under $\varphi_n$
of $x,y,z \in A$ respectively. Then the values $r_n$, $x_n$, $y_n$, and $z_n$
from Proposition \ref{finish} can be given by $1$, $x_n(\snum)$, 
$y_n(\snum)$, and $z_n(\snum)$ respectively.
Hence, if we set 
$a_n(\svar) = (y_n(\svar)+z_n(\svar))(\svar-1)$,
$b_n(\svar) = (x_n(\svar)+z_n(\svar))(\svar-1)$, and
$c_n(\svar) = (x_n(\svar)+y_n(\svar))(\svar-1)$, then 
it follows from Proposition \ref{finish} that both 1 and 2 of
Theorem \ref{maintheorem} are
satisfied for rational $\snum>1$. Note that if $\snum_0 \neq \snum_1$,
then $\Delta_n(\snum_0)$ is automatically not similar to
$\Delta_m(\snum_1)$ for any $m,n \geq 1$. To show that property 1
also holds for all {\em real} $\snum>1$, it suffices to show
$x_n(\snum), y_n(\snum), z_n(\snum) >0$. This follows from
continuity. 
\end{proofof}

\begin{corollary}
The set of rational points on $Y$ is Zariski dense in $Y$.
\end{corollary}
\begin{proof}
The infinitely many multiples of the section $R$ give infinitely many
curves on $Y$, each  with infinitely many rational points. Hence the
Zariski closure of the set of rational points is $Y$. 
\end{proof}

\begin{remark}
The four triples given in Remark \ref{concreteparam} 
correspond to the sections $R,3R,5R$, and $7R$. 
\end{remark}

\begin{remark}
As mentioned before, Randall Rathbun found with a
computer search a set of $8$ Heron triangles with the same area and
perimeter. His triangles correspond
to $r/(x+y+z)=28/195$ and the $8$ points on the corresponding elliptic curve
above $[28:195]$ generate a group of rank $4$. 
This yields relatively many points of
relatively low height. As in the proof of Proposition \ref{finish}
we can take any $n$ points on the real connected component that does not contain $\O$
and scale them to have the same perimeter and area. 
This is how we found the values in Table 1.
\check{$\tnum \leftarrow t$?}
\end{remark}

\bigskip

\mysection{\titlesix}{\label{secsix}}

In this section we will see how the N\'eron-Severi group of a
surface behaves under good reduction. Proposition \ref{nerred} is
known among specialists, but by lack of reference, we will 
include a proof, as sketched by B. Edixhoven. 
It will be used in the next
section to find the N\'eron-Severi group of $Y_{\overline{\Q}}$ and
the Mordell-Weil group $E(\overline{\Q}(\svar))$.

For the rest of this section, let $A$ be a discrete valuation ring of
a number field $L$ with residue field $k \isom \F_q$ for $q=p^r$. 
Let $S$ be an integral
scheme with a morphism $S \rightarrow \Spec A$ that is projective and 
smooth of relative dimension $2$. Then the 
projective surfaces $\overline{S} = S_{\overline{\Q}}$ and $\tilde{S} =
S_{\overline{k}}$ are smooth over the algebraically closed fields
$\overline{\Q}$ and $\overline{k}$ respectively.
We will assume that $\overline{S}$ and $\tilde{S}$ are integral, i.e.,
they are irreducible, nonsingular, projective surfaces. 

Let $l \neq p$ be a prime number. For any scheme $Z$ we set 
$$
H^i(Z_{\et},\Q_l) = \Q_l \otimes_{\Z_l} 
\left(\lim_{\leftarrow} H^i(Z_{\et},\Z/l^n\Z)\right).
$$
Furthermore, for every vector space $H$ over $\Q_l$ with a Galois
group $G(\overline{\F_q}/\F_q)$ 
acting on it, we define the twistings of $H$
to be the $G(\overline{\F_q}/\F_q)$-spaces $H(m)=H \otimes W^{\otimes m}$,
where 
$$
W=\Q_l \otimes_{\Z_l} (\lim_{\leftarrow} \mu_{l^n})
$$
is the one-dimensional $l$-adic vector space on which
$G(\overline{\F_q}/\F_q)$ operates according to its action on the group 
$\mu_{l^n}$ of $l^n$-th roots of unity. Here we use 
$W^{\otimes m} = \Hom(W^{\otimes -m}, \Q_l)$ for $m<0$.

For the rest of this article, all cohomology will be \'etale
cohomology, so we will often leave out the subscript $\et$.
In the rest of this section, we will use some results that are proven
in \myref{\milne} only in the noetherian case. 
Although we sometimes need the nonnoetherian case, besides
referring to the proofs of the general version, we will still also refer to
\myref{\milne}, because it is more easily available.

\begin{lemma}\label{cohombasechange}
Let $B$ denote the localization at some prime of the integral closure
of $A$ in $\overline{\Q}$.
For every integer $m$ the natural homomorphisms
$$
\eqalign{
H^i(S_B, \Q_l)(m) &\rightarrow H^i(\tilde{S}, \Q_l)(m) \qquad \mbox{and}\cr
H^i(S_B, \Q_l)(m) &\rightarrow H^i(\overline{S}, \Q_l)(m)
}
$$
are isomorphisms.
\end{lemma}
\begin{proof}
As tensoring with $W$ is exact, it suffices to prove this for $m=0$.
The ring $B$ is an integrally closed local ring
in an algebraically closed field with residue field isomorphic to
$\overline{k}$, so $B$ is a (nonnoetherian) strict Henselian ring 
(for the definition, see \myref{\egafour}, 
D\'ef. 18.8.2, or \myref{\milne}, \parag{} I.4).  
The surfaces $\overline{S}$ and $\tilde{S}$ are the fibers of 
$S_B \rightarrow \Spec B$. As $B$ is strictly Henselian, 
it follows from the proper base change theorem that the maps
$H^2(S_B,\Z/{l^n}\Z) \rightarrow H^2(\tilde{S},\Z/{l^n}\Z)$ are isomorphisms
for all $n\geq 0$,
see \myref{\milne}, Cor. VI.2.7, and \myref{\sgafourh}, p. 39,
Thm. IV.1.2. Hence, also the map $H^2(S_B,\Q_l) \rightarrow
H^2(\tilde{S}, \Q_l)$ obtained from taking the projective limit and
tensoring with $\Q_l$, is an
isomorphism. From the smooth base change theorem (\myref{\milne},
Thm. VI.4.1, and \myref{\sgafourh}, p. 63, Thm. V.3.2) it follows that 
$H^2(S_B,\Z/{l^n}\Z) \rightarrow H^2(\overline{S}, \Z/{l^n}\Z)$ is also an
isomorphism. For this exact statement, see \myref{\sgafourh},
p. 54--56: Lemme V.1.5, (1.6), and Variante. 
(for their $S$ take $S=\Spec B$; as $B$ is a strict Henselian local
ring which is integrally
closed in its algebraically closed fraction field already, we get
that their $S'$ equals their $S$). 
These statements assume that the morphism $S_B \rightarrow \Spec B$ is
locally acyclic, which follows from the fact that it is smooth, see 
\myref{\sgafourh}, p.58, Thm. (2.1).  
Passing to the limit and tensoring with $\Q_l$, we find that also the map 
$H^2(S_B,\Q_l) \rightarrow H^2(\overline{S}, \Q_l)$ is an isomorphism.
\donecheck{If I want to use ``any strict henselian ring in $\overline{\Q}$
dominating $A$'' for $B$, then I need that such a $B$ is a valuation
ring to apply valuative criterium and so that $B$ is integrally
closed, so that $S'=S$ in SGA 4.5, p.54-56. Then $\overline{S}$ is no
longer a fiber of $S_B \rightarrow \Spec B$, but that's ok.}
\end{proof}

\begin{proposition}\label{nerred}
There are natural injective homomorphisms 
\begin{equation}\label{nerredmaps}
\NS(\overline{S}) \otimes_{\Z_l} \Q_l \hookrightarrow
\NS(\tilde{S}) \otimes_{\Z_l} \Q_l \hookrightarrow
H^2(\tilde{S},\Q_l)(1)
\end{equation}
of finite dimensional vector spaces over $\Q_l$. The second injection
respects the Galois action of $G(\overline{k}/k)$. 
\end{proposition}
\begin{proof}
For any scheme $Z$, we have $H^1(Z_\et, \G_m) \isom \Pic Z$, see
\myref{\sgafourh}, p. 20, Prop. 2.3, or
\myref{\milne}, Prop. III.4.9. As long as $l \neq \charac k(z)$ for any
$z \in Z$, the Kummer sequence
$$
0 \rightarrow \mu_{l^n} \rightarrow \G_m \stackrel{[l^n]}{\longrightarrow}
\G_m \rightarrow 0 
$$
is a short exact sequence of sheaves on $Z_\et$, see 
\myref{\sgafourh}, p.21, (2.5), or \myref{\milne},
p. 66. Hence, from the long exact sequence we get a $\delta$-map
$$
\Pic Z \isom H^1(Z_\et, \G_m) \stackrel{\delta}{\longrightarrow} 
H^2(Z_\et, \mu_{l^n}).
$$
Taking the projective limit over $n$, this induces a homomorphism
$$
\Pic Z \rightarrow \lim_{\leftarrow} H^2(Z, \mu_{l^n}) \isom 
   \lim_{\leftarrow} H^2(Z, \Z/l^n\Z)\otimes \mu_{l^n} \rightarrow 
   H^2(Z, \Q_l)(1).
$$
Let $B$ be as in Lemma \ref{cohombasechange}.
Then from the above we get the diagram below, which commutes by
functoriality.  
Because $S_B$ is projective over $\Spec B$, it follows from the
valuative criterion for properness that the map $\Pic S_B \rightarrow
\Pic \overline{S}$ is an isomorphism (even if $f\colon X\rightarrow Y$
is a morphism between nonnoetherian schemes, it is still true that if
$f$ is proper, then for any valuation ring $A$ with quotient field $K$
the map $\Hom_Y(\Spec A, X)\rightarrow \Hom_Y(\Spec K, X)$ is a 
bijection, see \myref{\egaII}, Prop. 7.2.3 and Thm. 7.3.8). 
The bottom two maps are isomorphisms by Lemma \ref{cohombasechange}.

$$
\xymatrix{\Pic \overline{S} \ar[d]& \Pic S_B \ar[d]
               \ar[l]_{\isom} \ar[r] & \Pic \tilde{S} \ar[d] \\
  H^2(\overline{S},\Q_l)(1) &
 H^2(S_B,\Q_l)(1) \ar[l]_{\isom} \ar[r]^{\isom} &
 H^2(\tilde{S},\Q_l)(1) 
}
$$
%
%

For any proper variety $Z$ over an algebraically closed field, let 
$\Num Z$ denote the subgroup of $\Pic Z$ of all divisor classes on $Z$
that are numerically equivalent with $0$, i.e., those whose intersection 
number with every closed, integral curve on $Z$ is $0$, see 
\myref{\sgasix}, Exp. XIII, p. 644, 4.4. Then $\Pic Z/\Num Z$ is a
finitely generated free abelian group and in fact we have an isomorphism 
$\NS(Z)/\NS(Z)_{\tors} \isom \Pic Z/ \Num Z$, see
\myref{\hagtwo}, Prop. 3.1, and \myref{\sgasix}, Exp. XIII, p.645, Thm. 4.6. 
\donecheck{statement is for proper, alg. closed, 
defn. of num. equiv. should not be}
By \myref{\tatethree},
p. 97--98, the kernel of $\Pic Z \rightarrow H^2(Z, \Q_l)(1)$ is
$\Num Z$. From the diagram above, it follows that 
the composition
$$
\Pic \overline{S} \isom \Pic S_B \rightarrow \Pic \tilde{S} \rightarrow
H^2(\tilde{S},\Q_l)(1),
$$
which factors as
\begin{equation}\label{picToH2OfRed}
\eqalign{
\Pic \overline{S} \rightarrow \NS(\tilde{S})/\NS(\tilde{S})_{\tors} 
\hookrightarrow H^2(\tilde{S},\Q_l)(1)  \qquad \mbox{and as} \cr
\Pic \overline{S} \rightarrow
  H^2(\overline{S},\Q_l)(1) \isom
 H^2(S_B,\Q_l)(1) \isom  H^2(\tilde{S},\Q_l)(1), 
}
\end{equation}
has kernel $\Num \overline{S}$. Hence, the first map in
(\ref{picToH2OfRed}) induces an injective homomorphism
$\NS(\overline{S})/\NS(\overline{S})_{\tors}
\hookrightarrow\NS(\tilde{S})/\NS(\tilde{S})_{\tors}$. 
Tensoring with $\Q_l$ gives the desired homomorphisms.
\donecheck{Fulton Remark back? Then also his reference}
\end{proof}

%
%

For any variety $X$ over $k$
let $F_X\colon X \rightarrow X$ denote the {\em absolute 
Frobenius of $X$}, which acts as the identity on points, 
and by $f \mapsto f^p$ on the structure sheaf. 
Set $\varphi = F_{S_k}^r$  
and let $\varphi^*$ denote the automorphism on 
$H^2(\tilde{S},\Q_l)$ induced by $\varphi \times 1$ acting on 
$S_k \times_k \overline{k} \isom \tilde{S}$.

\begin{corollary}\label{rankNSbound}
The ranks of $\NS(\tilde{S})$ and 
$\NS(\overline{S})$ are bounded from above by the number of
eigenvalues $\lambda$ of $\varphi^*$ for which $\lambda/q$ is a root
of unity, counted with multiplicity. 
\end{corollary}
\begin{proof}
By Proposition \ref{nerred}, any upper bound for the rank of
$\NS(\tilde{S})$ is an upper bound for the rank of
$\NS(\overline{S})$. 
For any $k$-variety $X$, the absolute Frobenius $F_X$ acts as the 
identity on the site $X_{\et}$. 
Hence, if we set $\overline{X}= X \times_k \overline{k}$, then 
$F_{\overline{X}}=F_{X} \times F_{\overline{k}}$ acts as the identity
on $H^i(\overline{X}, \Q_l(m))$ for any $m$, see \myref{\tatethree},
\parag{} 3. Therefore, $F_X = F_X \times 1$ and $F_{\overline{k}} = 1
\times F_{\overline{k}}$ act as each other's inverses.

Let $\sigma\colon x \mapsto x^{q}$ denote the canonical
topological generator of $G(\overline{k}/k)$. Then $\sigma =
F_{\overline{k}}^r$ and as we have $\tilde{S} \isom S_k \times_k
\overline{k}$, we find $\varphi \times \sigma = 
F_{S_k}^r \times F_{\overline{k}}^r = F_{\tilde{S}}^r$.
By the above we find that  
the induced automorphisms $\varphi^{*(m)}$ and $\sigma^{*(m)}$ on
$H^2(\tilde{S},\Q_l)(m)$ act as each other's inverses for any $m$. 

As every divisor on $\tilde{S}$ is defined over a finite field
extension of $k$, some power of $\sigma^{*(1)}$ acts as the identity on
$\NS(\tilde{S})\subset H^2(\tilde{S},\Q_l)(1)$. 
It follows from Proposition \ref{nerred} that 
an upper bound for $\rk \NS(\tilde{S})$ is given by the number of
eigenvalues (with multiplicity) of $\sigma^{*(1)}$ 
that are roots of unity. 
As $\sigma^*$ acts as multiplication by $q$ on $W$, this equals the
number of eigenvalues $\nu$ of $\sigma^{*(0)}$ for which $\nu q$ is a
root of unity. The corollary follows as 
$\varphi^*=\varphi^{*(0)}$ acts as the inverse of $\sigma^{*(0)}$.
\end{proof}

\begin{remark}
Tate's conjecture states that the upper bound mentioned is
actually equal to the rank of $\NS(\tilde{S})$, see \myref{\tatethree}.
Tate's conjecture has been proven for ordinary K3 surfaces over fields
of characteristic $\geq 5$, see \myref{\nyog}, Thm. 0.2. 
In the next section we will be in exactly that situation.
\end{remark} 

\bigskip

\mysection{\titleseven}\label{secseven}

Note that also in this section all cohomology is \'etale cohomology,
so we often will leave out the subscript $\et$.
Let $\oY$ and $\oC$ denote
$Y_{\overline{\Q}}$ and $C_{\overline{\Q}}$ respectively, and set
$L=k(\oC) \isom \overline{\Q}(\svar) \supset \Q(\svar) = k(C) =K$.
By Theorem \ref{mordellrat} and Proposition \ref{finish} the points
$Q$ and $R$ both have infinite order in $E(L)$. 
Suppose there are integers $m,n$ such that $mQ+nR=0$.
Since complex conjugation sends $Q$ and $R$ to $-Q$ and $R$
respectively, we find that also $-mQ+nR=0$, whence $2mQ=2nR=0$. 
Therefore $m=n=0$, so $Q$ and $R$ are linearly independent, and $P$,
$Q$, and $R$ generate a group isomorphic to
$\Z^2 \times \Z/3\Z$. We will show that this is
the full Mordell-Weil group $E(L)$.  

\begin{proposition}\label{picard18}
The surface $\overline{Y}$ is a K3 surface. Its
N\'eron-Severi lattice has rank $\rho=18$.
The rank $r$ of the Mordell-Weil group $\overline{Y}(\overline{C})
\isom E(L)$ equals $r=2$.
\end{proposition}

\begin{proof}
To prove that $\overline{Y}$ is a K3 surface, it suffices by
definition to show that its
irregularity $q=H^1(\overline{Y},\O_{\overline{Y}})$ satisfies $q=0$ and
that any canonical divisor $K_{\oY}$ is linearly equivalent to $0$.

By Lemma \ref{samepiczero} we get $\Pic^0 \overline{Y}\isom \Pic^0
\overline{C}=0$, as $C$ is isomorphic to $\P^1$.
We conclude that 
$\NS(\overline{Y}) \isom \Pic(\overline{Y})$, so
algebraic and numerical
equivalence on $\overline{Y}$ coincide with linear equivalence. 
As $\tilde{X}$ is rational, we have $\chi(\O_{\tilde{X}})=
\chi(\O_{\P^2})=1$, see \myref{\hag}, Cor. V.5.6. By Proposition
\ref{baseextendellipticsurface} we get $\chi(\O_{\oY})=(\deg
\tilde{f}|_C) \cdot \chi(\O_{\tilde{X}_{\overline{\Q}}}) = 2$. From 
Theorem \ref{canonicalonminimalfibration} we then find that  
$K_{\oY} = 0$ in $\Pic \oY$. 
Hence, the canonical sheaf
$\omega_{\oY}$ is isomorphic to $\O_{\oY}$. We find from Serre
duality that $H^2(\oY,\O_{\oY}) \isom H^0(\oY, \omega_{\oY}) 
\isom H^0(\oY,\O_{\oY})$. Since $\oY$ is connected and projective,
we get $\dim H^2(\oY,\O_{\oY})=\dim H^0(\oY, \O_{\oY})=1$. 
Therefore, we get $q= \dim H^0(\oY,\O_{\oY})+ \dim
H^2(\oY, \O_{\oY}) - \chi(\O_{\oY})=1+1-2=0$. 

As seen in the proof of Proposition \ref{baseextendellipticsurface},
the singular fibers of $g$ come in pairs of copies of a singular fiber
of $\tilde{f}$. Hence, from Remark \ref{singfib} and Proposition
\ref{Tkern} we find that $\rho = 2 + 2 \left( (6-1) + (3-1) +
(1-1) + (1-1) \right) + r = 16+r$. Since $Q$ and $R$ are linearly
independent, we have $r \geq 2$, so we get $\rho \geq 18$.

We will show that $\rho \leq 18$ by reduction modulo a prime of good
reduction. Take $p=11$ and let $A=\Z_{(p)}$ be the localization of
$\Z$ at $p$ with residue field $k=A/p\isom \F_p$. 
Let $\X$ be the closed subscheme of $\P^3_A$ given by 
$r^2(x+y+z)=xyz$ and $\f\colon \X \dashrightarrow \P^1_A$ the
rational map that sends $[r:x:y:z]$ to $[r:x+y+z]$. 

As $\X$ is projective and $\X_\Q \isom X$, there are $A$-points $\Mi_i$ and
$\Ni_i$ on $\X$ such that $(\Ni_i)_\Q = N_i$ and $(\Mi_i)_\Q=M_i$. Let
$\pi' \colon
\tilde{\X} \rightarrow \X$ be the blow-up at the $6$ points $\Ni_i$ and
$\Mi_i$, and let $\tilde{\f} \colon \tilde{\X} \rightarrow \P^1_A$ be
the morphism induced by the composition $\f \circ \pi'$.
Let $\Ci \subset \tilde{\X}$ be the strict transform of the 
curve in $\X$ parametrized by 
$$
[r:x:y:z]=[\svar-1\mcol \svar+1\mcol \svar-1\mcol \svar(\svar-1)].
$$
Let $\Y$ denote the fibered product $\Y = \Ci \times_{\P^1_A}
\tilde{\X}$, and let $\g$ denote the projection $\Y\rightarrow \Ci$.
Then $\Y$ is a model of $Y$ over $A$, i.e., $\Y_{\Q} \isom Y$. 
Note that $\oY \isom \Y_{\overline{\Q}}$. Set 
$\tilde{Y}=\Y_{\overline{k}}$ and $\tilde{C}=\Ci_{\overline{k}}$. 
The following diagram shows how the base changes of $\Y$ that we will
deal with are related. A similar diagram holds for $\Ci$.


$$
\xymatrix{\overline{Y}\ar[d]_\isom&Y\ar[d]_\isom&&&\tilde{Y}\ar[d]_\isom\\
\Y_{\overline{Q}}\ar[r]\ar[d]&\Y_\Q\ar[r]\ar[d]&\Y\ar[d]&\Y_k\ar[l]\ar[d]&
\Y_{\overline{k}}\ar[l]\ar[d]\\  
\Spec\overline{\Q}\ar[r]&\Spec\Q\ar[r]&\Spec A&\Spec k\ar[l]&\Spec
\overline{k}\ar[l]\\ 
}
$$

We will show that $\Y$ is smooth over $\Spec A$.
Note that for each of the $\Ni_i$ and $\Mi_i$ there is an affine
neighborhood $U = \Spec S \subset \X$ for some $A$-algebra $S$, on which that
point corresponds to an ideal $I \subset S$ satisfying $pS \cap I^n =
pI^n$.  Set $T = S \otimes_A k \isom  S/pS$ and $J = IT$.
Then $U_k = \Spec T$ and we have 
$$
I^n \otimes_A k \isom I^n/pI^n \isom I^n/(pS \cap I^n) \isom 
I^n \cdot S/pS \isom I^nT = J^n.
$$
This implies 
$$
\Proj \left( T \oplus J \oplus J^2 \oplus \ldots \right) \isom
\Proj\left(S \oplus I \oplus I^2 \oplus \ldots\right) \times_{\Spec A} \Spec k,
$$
which tells us that the blow-up of the reduction $\X_k$
at the points $(\Mi_i)_k$ and $(\Ni_i)_k$
is isomorphic to $\tilde{\X}\times_A
k$, i.e., the reduction $\tilde{\X}_k$ of $\tilde{\X}$. 

One easily checks that $\X_k$ is geometrically regular
outside the three ordinary double points $(\Mi_i)_k$. Hence, this
blow-up of $\X_k$ at the points $(\Mi_i)_k$ and $(\Ni_i)_k$ is smooth
over $k$, see \myref{\hag}, exc. I.5.7. 
Thus $\tilde{\X}_k$ is smooth over $k$. 
As the morphism $\Ci_k \rightarrow \P^1_k$ is unramified at
the points of $\P^1_k$ where $\tilde{\f}_k$ has singular fibers 
(as is easily checked), $\Y_k$ is smooth over $k$ as well
(cf. Prop. \ref{baseextendellipticsurface}). Since the other
fiber $\Y_{\Q}\isom Y$ of $\Y \rightarrow \Spec A$ is also smooth 
over its ground field $\Q$, we
conclude that $\Y$ is smooth over $\Spec A$ (cf. Remark \ref{smooth}).

Let $\varphi \colon \Y_k \rightarrow \Y_k$ denote the
absolute Frobenius of $\Y_k$ as in the previous section. 
Let $\varphi_i^*$ denote the induced
automorphism on $H^i(\tilde{Y}, \Q_l)$. By Corollary \ref{rankNSbound}
the Picard number $\rho$ is bounded from above by the number of
eigenvalues $\lambda$ of $\varphi_2^*$ for which $\lambda/p$ is a root
of unity. 
We will count these eigenvalues using 
the Lefschetz trace formula and the Weil conjectures. 
The characteristic polynomial of $(\varphi^*_i)^n$ acting on
$H^i(\tilde{Y},\Q_l)$ is 
$$
P_i(t) = \det \left( \Id - t \cdot (\varphi^*_i)^n|H^i(\tilde{Y},\Q_l) 
\right) = \prod_{i=1}^{b_i} (1-\alpha_{ij}t).
$$
By the Weil conjectures, $P_i(t)$ is a rational polynomial and the
reciprocal roots have absolute value $|\alpha_{ij}|=p^{ni/2}$, see
\myref{\deligne}, Thm. 1.6. 


By Lemma \ref{cohombasechange} we have 
$\dim H^i(\oY,\Q_l) = \dim H^i(\tilde{Y},\Q_l)$ for $0 \leq i \leq 4$.
Since $\oY$ is a K3 surface, the betti numbers equal
$\dim H^i(\tilde{Y},\Q_l)=b_i=1,0,22,0,1$ for $i=0,1,2,3,4$
respectively. Therefore, from the Weil conjectures we find 
$P_i(t)=1-t,1,1,1-p^{2}t$ for \donecheck{why not $1+t$ or $1+p^2t$?}
$i=0,1,3,4$ respectively, whence $\Tr \varphi^*_i = 1,0,0,p^2$ for
$i=0,1,3,4$. Similarly, we get $\Tr (\varphi^*_i)^n = 1,0,0,p^{2n}$ for
$i=0,1,3,4$ and $n\geq 1$. That means that for any $n\geq 1$, if we
know the number of $\F_{p^n}$-points of $\Y_k$, then from the
Lefschetz Trace Formula (see \myref{\milne}, Thm. VI.12.3).
$$
\# \Y_k(\F_{p^n}) = \sum_{i=0}^4 (-1)^i 
\Tr\left((\varphi^*_i)^n|H^i(\tilde{Y},\Q_l)\right)
$$
we can compute $\Tr (\varphi^*_2)^n|H^2(\tilde{Y},\Q_l)$.

Let $V$ denote the image in $H^2(\tilde{Y},\Q_l)$ under the composition in
(\ref{nerredmaps}) of the $18$-dimensional subspace of 
$\NS(\oY) \otimes \Q_l$
that we already know, i.e., generated by the irreducible components of
the singular fibers of $g$ and the sections $\O$, $Q$, and $R$. 

All these generators of $V$ are defined over the $k=\F_p$, except the image of
$Q$, which is defined over $\F_{p^2}$. In the Mordell-Weil group
modulo torsion $\tilde{Y}(\tilde{C})/\tilde{Y}(\tilde{C})_\tors$ we
have $\varphi (Q) = -Q$. Hence $V$ is $\varphi^*_2$-invariant and we 
find that $\Tr (\varphi_2^*)^n|V = 17p^n+(-1)^np^n$. 

Let $W$ be any subspace of 
$H^2(\tilde{Y},\Q_l)$, such that $V \oplus W \isom
H^2(\tilde{Y},\Q_l)$, and let
$\pi_W$ denote the projection of $H^2(\tilde{Y},\Q_l)$ onto $W$. 
Then $W$ has dimension $4$ and from just linear algebra we get
\begin{equation}\label{charpoly}
\charpoly (\varphi^*_2|H^2) = \charpoly(\varphi^*_2|V) \cdot 
\charpoly(\pi_W\circ\varphi_2^*|W)
\end{equation}
and
$$
\Tr ((\varphi^*_2)^n|H^2) = \Tr((\varphi^*_2)^n|V) +
\Tr(\pi_W\circ(\varphi^*_2)^n|W)
$$
The last equality then allows us to compute
$\Tr(\pi_W\circ(\varphi^*_2)^n|W)$ for $n\geq 
1$, which is done for $n=1,2,3$ in the following table.

\begin{center}
$$
\begin{array}{|l||r|r|r|}
\hline
n                                         & 1 & 2 & 3 \cr
\hline
\Tr (\varphi^*_0)^n|H^0(\tilde{Y},\Q_l) & 1 & 1 & 1 \cr
\Tr (\varphi^*_1)^n|H^1(\tilde{Y},\Q_l) & 0 & 0 & 0 \cr
\Tr (\varphi^*_3)^n|H^3(\tilde{Y},\Q_l) & 0 & 0 & 0 \cr
\Tr (\varphi^*_4)^n|H^4(\tilde{Y},\Q_l) & p^2 & p^4 & p^6 \cr
\# \Y_k(\F_{p^n})                       & 298 & 16908 & 1792858 \cr
\Tr (\varphi^*_2)^n|H^2(\tilde{Y},\Q_l) & 176 & 2266 & 21296 \cr
\Tr (\varphi^*_2)^n|V                   & 16p & 18p^2 & 16p^3 \cr
\Tr (\pi_W \circ (\varphi^*_2)^n|W)     & 0 & 88 & 0 \cr
\hline
\end{array}
$$
\smallskip

Table 2
\end{center}

We computed the number of points on $\Y_k(\F_{p^n})$ as follows.
As $\Y_k$ has the structure of elliptic surface
over $\Ci_k$, we can let the computer package Magma compute the 
number of points above every point of 
$\Ci_k(\F_{p^n})$ with a nonsingular elliptic fiber. 
Adding to that the contribution of the
singular fibers gives the total number of points.

%

For any linear operator $T$ on an $m$-dimensional vector space with
characteristic polynomial 
$$
\charpoly T = X^m+c_1 X^{m-1} +c_2 X^{m-2}+\ldots+c_{m-1}X+c_m,
$$
we have $c_1 = -t_1$, $c_2 = \frac{1}{2}(t_1^2-t_2)$, and 
$c_3 = -\frac{1}{6}(t_1^3+2t_3-3t_1t_2)$, where $t_n = \Tr T^n$. 
From this and table 2 we find that the characteristic
polynomial of $\pi_W \circ \varphi^*_2 |W$ equals $h=X^4-44X^2+c_4$ for
some $c_4$. By the Weil conjectures, and (\ref{charpoly}), 
the roots of $h$ have absolute 
value $p$ and their product $c_4$ is rational, so $c_4=\pm p^4$. 
As not all roots of $X^4-44X^2-11^4$ have absolute value $11$, we get
$h=X^4-44X^2+11^4$. If $\alpha$ is a root of $h$ then
$\beta = (\alpha/p)^2$ satisfies $11\beta^2-4\beta+11=0$. 
As the only quadratic roots of unity are 
$\pm \sqrt{-1}$ and $\zeta_6^i$, we find that $\beta$, whence
$\alpha/p$, is not a root of unity. 
From (\ref{charpoly}) it follows that 
$\alpha/p$ is a root of unity for at most $22-4=18$
roots $\alpha$ of $\charpoly (\varphi^*_2|H_2)$. From Corollary
\ref{rankNSbound} we find $\rho \leq 18$.
\end{proof}

\begin{corollary}\label{MW}
The Mordell-Weil group $E(L)$ is generated by 
$P$, $Q$, and $R$ and is isomorphic to $\Z^2 \times \Z/3\Z$.
\end{corollary}
\begin{proof}
As $\overline{Y}\rightarrow \overline{C}$ is a relatively minimal
fibration and $\overline{Y}$ is regular and projective, the 
N\'eron model of $\overline{Y}/\overline{C}$ is  
obtained from $\overline{Y}$ by
deleting the singular points of the singular fibers, see
\myref{\silvtwo}, Thm. IV.6.1, and \myref{\neron}, 
\parag{} 1.5, Prop. 1. 
Note that at $\svar=0$ and $\svar+1=0$ we have additive reduction (type
IV), whence the identity component of the reduction has no torsion.
Since we are in characteristic $0$,
the kernel of reduction $E_1(L)$ has no torsion
either, see \myref{\silv}, Prop. VII.3.1. It follows that the group
$E_0(L)$ of nonsingular reduction has no torsion, see \myref{\silvtwo},
Rem. IV.9.2.2. By the classification of singular fibers we find that 
$E(L)/E_0(L)$ has order at most $3$, see \myref{\silvtwo},
Cor. IV.9.2 and Tate's Algorithm IV.9.4. We
conclude that $E(L)_{\tors}$ has order $3$ and is generated by $P$.
%

With Shioda's explicit formula for the Mordell-Weil pairing
(\myref{\shioda}, Thm. 8.6), we find that $\langle Q,R \rangle =0$ and 
$\langle Q,Q \rangle =\langle R,R \rangle =1$. Hence, as seen before,
$Q$ and $R$ are
linearly independent. As the rank $r=\rk E(L)$ equals $2$ by
Proposition \ref{picard18}, the group generated by $Q$ and $R$ has
finite index in the Mordell-Weil lattice $E(L)/E(L)_{\tors}$.
If the Mordell-Weil lattice were not generated by $Q$ and $R$, then
it would contain a nonzero element $S=aQ+bR$ with $a,b \in \Q$ and
$-\frac{1}{2}<a,b \leq \frac{1}{2}$, so that $\langle S, S \rangle = 
a^2+b^2 \leq \frac{1}{2}$. On the other hand, based on the type of
singularities, it follows from the
explicit formulas for the Mordell-Weil pairing that its values 
are contained in $\frac{1}{6}\Z$. As for any rational $a,b$ the
$3$-adic valuation of $a^2+b^2$ is even, we conclude that in fact we
have $\langle S, S \rangle \in \frac{1}{2}\Z$, so that $a^2+b^2 \geq
\frac{1}{2}$. Thus, we find $a^2+b^2=\frac{1}{2}$, whence
$a=b=\frac{1}{2}$. Therefore, $2S =
Q+R+\varepsilon P$ for some $\varepsilon\in\{0,1,2\}$. After adding
$\varepsilon P$ to $S$ if necessary, we may assume $\varepsilon =0$ 
without loss of generality. \donecheck{or use discriminants to show that
  the index is at most $2$. (probably barely faster)}

It suffices to check that $Q+R \not \in 2E(K)$. Let $(p_S,q_S)$,
$(p_{2S},q_{2S})$, and $(p_{Q+R},q_{Q+R})$ denote
the Weierstrass coordinates of $S$, $2S$, and $Q+R$ respectively. 
Using addition formulas, we can compute 
$p_{Q+R}\in \Q(i)(\svar)$ explicitely and express $p_{2S}$ in terms of $p_S$. 
Let $u$ be defined by $p_S-4(\svar-1)^2=2(\svar-1)u$. 
Then in terms of $u$, the equation $p_{2S} = p_{Q+R}$
simplifies to 
\begin{equation}\label{pofS}
\eqalign{
u^4& + 4(\svar-1)(\svar+1)(\svar+i)u^3 + 
               2(\svar^2+(1+i)\svar-2+i)\svar^2(\svar+1)^2u^2 + \cr
      &8(\svar^2+(1+i)\svar-2+i)(\svar-1)\svar^2(\svar+1)^2u + 
               8(\svar+i)\svar^2(\svar-1)^2(\svar+1)^3 =0\cr
}
\end{equation}
By Gauss's Theorem any root $u\in L=\overline{\Q}(\svar)$ of this
equation is contained in $\overline{\Q}[\svar]$ and divides the
constant term $8(\svar+i)\svar^2(\svar-1)^2(\svar+1)^3$. Hence, any root
$u$ is of the form 
$$
u = c \svar^k (\svar+1)^l (\svar-1)^m (\svar+i)^n,
$$
for some constant $c$ and exponents $k$, $l$, $m$, and $n$.
Considering the four Newton polygons, we find 
$k=0$, $l=1$, and $m,n \in \{0,1\}$.
One easily checks that for none of the four possibilities for $m,n$
there is a $c$ such that (\ref{pofS}) is satisfied.
\end{proof}

\begin{corollary}
The discriminant of the N\'eron-Severi lattice 
$\NS(\oY)$ equals $-36$.
\end{corollary}
\begin{proof}
From the short exact sequence (\ref{TNsE}) we find the following
equation, relating the discriminant of the N\'eron-Severi lattice to
that of the Mordell-Weil lattice, see \myref{\shioda}, Thm. 8.7.
$$
|\disc \NS(\oY)| = \frac{\disc MW \cdot \prod m_v^{(1)}}{|E(L)_\tors|^2},
$$  
where $m_v^{(1)}$ is the number of irreducible components of
multiplicity $1$ of the fiber of $g$ above $v\in C$, 
and $MW$ stands for the Mordell-Weil lattice
$E(L)/E(L)_\tors$. Note that we used
$\disc T = \prod m_v^{(1)}$, see \myref{\shioda}, (7.9). 
In the proof of Corollary \ref{MW} we
have seen that $\disc MW = 1$, so we get 
$$
|\disc \NS(\oY)| = \frac{1 \cdot 6 \cdot 6 \cdot 3 \cdot 3}{3^2} = 36.
$$
By the Hodge index Theorem $\disc \NS(\oY)$ is negative, so we
get $\disc \NS(\oY) = -36$. \donecheck{could conclude that $\overline{Y}$
  is not Kummer, but doesn't really seem worth it.}
\end{proof}

\bigskip

\mysectiononumber{\titlenine}\label{secnine}

\luijkreflist
\myrefart{\aas}{Aassila, M.}{Some results on Heron
  triangles}{Elem. Math.}{{\bf 56} (2001)}{143--146}
\myrefart{\artin}{Artin, M.}{On Isolated Rational Singularities of
  Surfaces}{Amer. J. Math.}{{\bf 88} (1966)}{129--136}
\myrefart{\bogomolovtschink}{Bogomolov, F. and Tschinkel, Yu.}{Density of rational
  points on elliptic K3 surfaces}{Asian J. Math.}{{\bf 4}, 2
  (2000)}{351--368}  
\myrefbook{\neron}{Bosch, S., L\"utkebohmert, W., and Raynaud,
  M.}{N\'eron Models}{Springer-Verlag, Berlin}{1990}
\myrefart{\bombmum}{Bombieri, E. and Mumford, D.}{Enriques' classification of
surfaces in char. $p$, II}{Complex Analysis and Algebraic
Geometry--Collection of papers dedicated to K. Kodaira}{ed. W.L. Baily and T.
Shioda, Iwanami and Cambridge Univ. Press (1977)}{23--42}
\myrefart{\bruce}{Bruce, J. and Wall, C.}{On the classification of
  cubic surfaces}{J. London Math. Soc. (2)}{{\bf 19} (1979)}{245--256}
\myrefart{\chin}{Chinburg, T.}{Minimal Models of Curves over Dedekind
rings}{Arithmetic Geometry}{ed. Cornell, G. \& Silverman,
  J. (1986)}{309--326} 
\myrefart{\deligne}{Deligne, P.}{La Conjecture de
  Weil. I}{Publ. Math. IHES}{{\bf 43} (1974)}{273--307}
\myrefart{\duval}{Du Val, P.}{On isolated singularities which do not
  affect the conditions of adjunction, Part I}{Proc. Cambridge
  Phil. Soc.}{{\bf 30} (1934)}{453--465}
\myrefbook{\egaII}{Grothendieck, A.}{\'El\'ements de g\'eom\'etrie
  alg\'ebrique. IV. \'Etude globale
  \'el\'ementaire de quelques classes de morphismes}
{Inst. Hautes \'Etudes Sci. Publ. Math., no. {\bf 8}}{1961}
\myrefbook{\egaone}{Grothendieck, A.}{\'El\'ements de g\'eom\'etrie
  alg\'ebrique. IV. 
\'Etude locale des sch\'emas et des morphismes de sch\'emas, Premi\`ere
  partie}{Inst. Hautes 
\'Etudes Sci. Publ. Math., no. {\bf 20}}{1964}
\myrefbook{\egatwo}{Grothendieck, A.}{\'El\'ements de g\'eom\'etrie
  alg\'ebrique. IV.
\'Etude locale des sch\'emas et des morphismes de sch\'emas, Seconde
  partie}{Inst. Hautes 
\'Etudes Sci. Publ. Math., no. {\bf 24}}{1965}
\myrefbook{\egafour}{Grothendieck, A.}{\'El\'ements de g\'eom\'etrie
  alg\'ebrique. IV.
\'Etude locale des sch\'emas et des morphismes de sch\'emas, Quatri\`eme
  partie}{Inst. Hautes 
\'Etudes Sci. Publ. Math., no. {\bf 32}}{1967}
\myrefart{\hagtwo}{Hartshorne, R.}{Equivalence relations of algebraic
  cycles and subvarieties of small codimension}{Algebraic Geometry,
  Arcata 1974}{Amer. Math. Soc. Proc. Symp. Pure Math. {\bf
  29} (1975)}{129--164} 
\myrefbook{\hag}{Hartshorne, R.}{Algebraic Geometry}{GTM {\bf 52},
Springer-Verlag, New-York}{1977}
\myrefart{\kramluc}{Kramer, A.-V. and Luca, F.}{Some remarks on Heron
  triangles}{Acta Acad. Paedagog. Agriensis Sect. Mat. (N.S.)}{{\bf
  27}  (2000)}{25--38 (2001)} 
\myrefart{\kodtwo}{Kodaira, K.}{On compact analytic surfaces II-III}{Ann. of
Math.}{{\bf 77} (1963), pp. 563--626; {\bf 78} (1963)}{1--40}
\myrefart{\kod}{Kodaira, K.}{On the structure of compact complex
  analytic surfaces I, II}{Amer. J. Math.}{{\bf 86} (1964),
  pp. 751--798; {\bf 88} (1966)}{682--721} 
\myrefart{\lich}{Lichtenbaum, S.}{Curves over discrete valuation
  rings}{Amer. J. Math.}{{\bf 90} (1968)}{380--405}
\myrefart{\lip}{Lipman, J.}{Rational singularities, with applications
  to algebraic surfaces and unique
  factorization}{IHES Publ. Math.}{{\bf 36} (1969)}{195--279}
\myrefbook{\manin}{Manin, Y.}{Cubic forms: Algebra, Geometry,
  Arithmetic}{North-Holland, Amsterdam}{1974}
\myrefart{\maz}{Mazur, B.}{Modular curves and the Eisenstein ideal}{IHES Publ.
Math.}{{\bf 47} (1977)}{33--186}
\myrefbook{\milne}{Milne, J.S.}{\'Etale Cohomology}{Princeton
  Mathematical Series 
  {\bf 33}, Princeton University Press, New Jersey}{1980}
\myrefart{\nagata}{Nagata, M.}{On rational surfaces I,
  II}{Mem. coll. Sci. Kyoto (A)}{{\bf 32} (1960), pp. 351--370; {\bf
    33} (1960)}{271--293}
\myrefart{\nyog}{Nygaard, N. and Ogus, A.}{Tate's conjecture for K3
  surfaces of finite height}{Annals of Math.}{{\bf 122} (1985)}{461--507}
\myrefart{\pink}{Pinkham, H.}{Singularit\'es Rationelles de
  Surfaces}{S\'eminaire sur les Singularit\'es des Surfaces}
{Lect. Notes in Math. {\bf 777},
ed. M. Demazure, H. Pinkham, and B. Teissier, Springer-Verlag
  (1980)}{147--172} 
\myrefbook{\sull}{O'Sullivan, M.}{Classification and Divisor Class Groups of
  Normal Cubic Surfaces in $\P^3$}{U.C. Berkeley, Ph.D. dissertation
  (not published)}{1996}
\myrefbook{\sgaone}{Grothendieck, A.}{Rev\^etements \'etales et Groupe
Fondamental}{Lect. Notes in Math. {\bf 224}, Springer-Verlag,
  Heidelberg}{1971} 
\myrefbook{\sgafourh}{Grothendieck, A. et al.}{Cohomologie \'etale}
{Lect. Notes in Math. {\bf 569}, Springer-Verlag,
  Heidelberg}{1977} 
\myrefbook{\sgasix}{Grothendieck, A. et al.}{Th\'eorie des
  Intersections et Th\'eor\`eme de Riemann-Roch}{Lect. Notes in
  Math. {\bf 225}, Springer-Verlag, Heidelberg}{1971} 
\myrefbook{\shaf}{Shafarevich, I.}{Lectures on Minimal Models and Birational
Transformations of Two-dimensional Schemes}{Tata Institute, Bombay}{1966}
\myrefart{\shioda}{Shioda, T.}{On the Mordell-Weil Lattices}{Comm. Math. Univ.
Sancti Pauli}{{\bf 39}, 2 (1990)}{211--240}
\myrefbook{\silv}{Silverman, J.H.}{The Arithmetic of Elliptic Curves}{GTM {\bf
106}, Springer-Verlag, New-York}{1986}
\myrefbook{\silvtwo}{Silverman, J.H.}{Advanced Topics in the Arithmetic of
Elliptic Curves}{GTM {\bf 151}, Springer-Verlag, New-York}{1994}
\myrefart{\tatetwo}{Tate, J.}{Genus change in inseparable extensions
  of function fields}{Proc. AMS}{{\bf 3} (1952)}{400--406}
\myrefart{\tatethree}{Tate, J.}{Algebraic cycles and poles of zeta
  functions}{Arithmetical Algebraic Geometry}
  {ed. O.F.G. Schilling (1965)}{93--110}
\myrefart{\tate}{Tate, J.}{Algorithm for determining the type of a 
singular fiber in an elliptic pencil}
{Modular functions of one variable IV}{Lect. Notes in Math. {\bf 476},
ed. B.J.~Birch and W.~Kuyk, Springer-Verlag, Berlin (1975)}{33--52}
\eindluijkreflist

\end{document}